\documentclass[12pt]{amsart}
\usepackage{amscd}      
\usepackage{amssymb}
\usepackage{amsmath}
\usepackage{xypic}      
\LaTeXdiagrams          
\usepackage[all,v2]{xy}
\xyoption{2cell} \UseAllTwocells \xyoption{frame} \CompileMatrices
\allowdisplaybreaks[3]

\usepackage{latexsym}
\usepackage{epsfig}
\usepackage{amsfonts}
\usepackage{enumerate}

\newtheorem{prop}{Proposition}[section]
\newtheorem{lem}[prop]{Lemma}

\newtheorem{cor}[prop]{Corollary}
\newtheorem{thm}[prop]{Theorem}

\numberwithin{equation}{section}

\newtheorem{defn}[prop]{Definition}

\newtheorem{example}[prop]{Example}

\newtheorem{rmk}[prop]{Remark}

\newenvironment{pf}{\begin{trivlist}\item[]{\sc Proof.}}%
            {\nolinebreak $\Box$ \end{trivlist}}

\newcommand{\noprint}[1]{}

\renewcommand{\tilde}{\widetilde}

\newcommand{\toto}{\rightrightarrows}

\newcommand{\ldiag}[1]%
       {\makebox[0cm]{${\scriptstyle#1}\downarrow\phantom{\scriptstyle#1}$}}
\newcommand{\ldiagup}[1]%
       {\makebox[0cm]{${\scriptstyle#1}\uparrow\phantom{\scriptstyle#1}$}}
\newcommand{\rdiag}[1]%
       {\makebox[0cm]{$\phantom{\scriptstyle#1}\downarrow{\scriptstyle#1}$}}
\newcommand{\sediagr}[1]%
       {\makebox[0cm]{$\phantom{\scriptstyle#1}\searrow{\scriptstyle#1}$}}
\newcommand{\nediagr}[1]%
       {\makebox[0cm]{$\phantom{\scriptstyle#1}\nearrow{\scriptstyle#1}$}}
\newcommand{\rdiagup}[1]%
       {\makebox[0cm]{$\phantom{\scriptstyle#1}\uparrow{\scriptstyle#1}$}}
\newcommand{\swdiag}[1]%
       {\makebox[0cm]{$\phantom{\scriptstyle#1}\swarrow{\scriptstyle#1}$}}
\newcommand{\sediag}[1]%
       {\makebox[0cm]{${\scriptstyle#1}\searrow\phantom{\scriptstyle#1}$}}
\newcommand{\nediag}[1]%
       {\makebox[0cm]{${\scriptstyle#1}\nearrow\phantom{\scriptstyle#1}$}}

\newcommand{\doublearrowstack}[2]%
                      {{{{\scriptstyle#1}\atop{\textstyle\longrightarrow}}\atop{{\textstyle\longrightarrow}\atop{\scriptstyle#2}}}}
\newcommand{\rightleftarrowstack}[2]%
                      {{{{\scriptstyle#1}\atop{\textstyle\longrightarrow}}\atop{{\textstyle\longleftarrow}\atop{\scriptstyle#2}}}}
\newcommand{\leftrightarrowstack}[2]%
                      {{{{\scriptstyle#1}\atop{\textstyle\longleftarrow}}\atop{{\textstyle\longrightarrow}\atop{\scriptstyle#2}}}}

\newcommand{\overtoparrow}%
{\makebox[0cm]{\beginpicture \setcoordinatesystem units
<.8cm,.4cm> point at 0 0 \setplotarea x from -3 to 3, y from 0 to
1 \setquadratic \plot -3 0 0 1 3 0 / \put{\vector(3,-1){0}}[Bl] at
3 0
\endpicture}}

\newcommand{\underbottomarrow}%
{\makebox[0cm]{\beginpicture \setcoordinatesystem units
<.8cm,.4cm> point at 0 0 \setplotarea x from -3 to 3, y from 0 to
1 \setquadratic \plot -3 1 0 0 3 1 / \put{\vector(3,1){0}}[Bl] at
3 1
\endpicture}}

\newcommand{\ses}[5]%
{0\longrightarrow#1\stackrel{#2}{ \longrightarrow}#3\stackrel{#4}{
\longrightarrow}#5\longrightarrow0}

\newcommand{\dt}[6]%
{#1\stackrel{#2}{longrightarrow}#3
\stackrel{#4}{\longrightarrow}#5 \stackrel{#6}{\longrightarrow}
#1[1]}

\newcommand{\cat}[1]%
{(\mbox{\rm #1})}

\def\Label#1{\label{#1}{\tt [#1]}\phantom{h}}
\def\Label{\label}

\newcommand{\sX}{\mathcal{X}}

\newcommand{\sL}{\mathcal{L}}

\def\<{\left\langle}
\def\>{\right\rangle}

\title[Integral Chow rings for toric DM stacks]{The Integral (orbifold) Chow Ring of Toric Deligne-Mumford Stacks}

\author{Yunfeng Jiang}
\address{Department of Mathematics\\ University of British Columbia\\ 1984 Mathematics Road\\Vancouver\\ BC V6T 1Z2\\ Canada}
\email{jiangyf@math.ubc.ca}
\address{$Current~ Address$: Department of Mathematics\\ University of Utah\\ 155 S 1400 E JWB 233\\Salt Lake city\\ UT 84112\\ USA}
\email{jiangyf@math.utah.edu}
\author{Hsian-Hua Tseng}
\address{Department of Mathematics\\ University of British Columbia\\ 1984 Mathematics Road\\Vancouver\\ BC V6T 1Z2\\ Canada}
\email{hhtseng@math.ubc.ca}
\date{\today}

\begin{document}
\begin{abstract}
In this paper we study the integral Chow ring of toric
Deligne-Mumford stacks. We prove that the integral Chow ring of a
semi-projective toric Deligne-Mumford stack is isomorphic to the
Stanley-Reisner ring of the associated stacky fan. The integral
orbifold Chow ring is also computed. Our results are illustrated with several examples. 

\end{abstract}

\maketitle

\section{Introduction}

Chow groups with integer coefficients of algebraic stacks were defined by Edidin and Graham \cite{EG}
and Kresch \cite{Kr}, using Totaro's idea \cite{To} of integral Chow ring of classifying spaces. In \cite{EG}, the authors constructed an intersection theory of stack quotient $[X/G]$ of a quasi-projective variety $X$ by an algebraic group $G$. In the case of Deligne-Mumford stacks, the authors proved that the equivariant Chow ring $A^{*}_{G}(X)$ with integer coefficients is isomorphic to the integral Chow ring of the quotient stack $[X/G]$. 

Toric Deligne-Mumford stacks were introduced by Borisov, Chen and Smith \cite{BCS} via a generalization of the quotient construction \cite{Cox} of simplicial toric varieties. A construction of toric stacks using logarithmic geometry can be found in \cite{Iwa1}. The purpose of this paper is to compute the (orbifold) Chow ring with integer coefficients of a toric Deligne-Mumford stack. A toric Deligne-Mumford stack is defined in terms of a stacky fan
$\mathbf{\Sigma}=(N,\Sigma,\beta)$, where $N$ is a finitely
generated abelian group, $\Sigma\subset N_\mathbb{Q}=N\otimes_{\mathbb{Z}}\mathbb{Q}$ is a simplicial fan
and $\beta: \mathbb{Z}^{n}\to N$ is a map determined
by the elements $\{b_{1},\cdots,b_{n}\}$ in $N$.  By assumption, $\beta$ has finite cokernel and
$\{\overline{b}_{1},\cdots,\overline{b}_{n}\}$ generate the
simplicial fan $\Sigma$, where $\overline{b}_{i}$ is the image of
$b_{i}$ under the natural map $N\to N_\mathbb{Q}$.
The  toric Deligne-Mumford stack $\mathcal{X}(\mathbf{\Sigma})$  associated to $\mathbf{\Sigma}$
is defined to be  the quotient stack
$[Z/G]$, where $Z$ is the open subvariety $\mathbb{C}^{n}\setminus\mathbb{V}(J_{\Sigma})$, $J_{\Sigma}$
is the irrelevant ideal of the fan, $G$ is the product of an algebraic torus and a finite abelian group.
The $G$-action on $Z$ is given via a group homomorphism $\alpha: G \to (\mathbb{C}^{*})^{n}$, where
$\alpha$ is obtained by taking $\text{Hom}_{\mathbb{Z}}(-,\mathbb{C}^{*})$ functor to the Gale dual
$\beta^{\vee}: \mathbb{Z}^{n}\to N^{\vee}$ of $\beta$ and $G=\text{Hom}_{\mathbb{Z}}(N^{\vee},\mathbb{C}^{*})$.

Each ray $\rho_{i}$ in the fan $\Sigma$ gives a line bundle $\mathcal{L}_{i}$
over $\mathcal{X}(\mathbf{\Sigma})$ which is defined by the quotient
$Z\times \mathbb{C}/G$ and the action of $G$ on $\mathbb{C}$ is through the
$i$-th component of the map $\alpha$. The Picard group of $\mathcal{X}(\mathbf{\Sigma})$ is seen to be isomorphic to $N^{\vee}$.

Every stacky fan $\mathbf{\Sigma}$ has an underlying {\em reduced} stacky
fan $\mathbf{\Sigma_{red}}=(\overline{N},\Sigma,\overline{\beta})$,
where $\overline{N}=N/\text{torsion}$, $\overline{\beta}:
\mathbb{Z}^{n}\to \overline{N}$ is the natural projection given by the vectors
$\{\overline{b}_{1},\cdots,\overline{b}_n\}\subseteq \overline{N}$.
The toric Deligne-Mumford stack $\mathcal{X}(\mathbf{\Sigma_{red}})$
is a toric orbifold. By construction 
$\mathcal{X}(\mathbf{\Sigma_{red}})=[Z/\overline{G}]$, where
$\overline{G}=\text{Hom}_{\mathbb{Z}}(\overline{N}^{\vee},\mathbb{C}^{*})$
and $\overline{N}^{\vee}$ is the Gale dual 
$\overline{\beta}^{\vee}: \mathbb{Z}^{n}\to \overline{N}^{\vee}$
of the map $\overline{\beta}$. The stack $\mathcal{X}(\mathbf{\Sigma_{red}})$ can be obtained by the rigidification construction (see e.g. \cite{ACV}). 
Each ray $\rho_{i}$ in the fan $\Sigma$ also gives a line bundle $L_{i}$
over the toric orbifold which corresponds to the divisor $D_{i}$. This line bundle
$L_{i}$ also has a quotient construction $Z\times \mathbb{C}/\overline{G}$ where $\overline{G}$ acts on $\mathbb{C}$ via the $i$-th component of the map $\overline{\alpha}: \overline{G}\to (\mathbb{C}^{*})^{n}$, obtained by taking $\text{Hom}_\mathbb{Z}(-,\mathbb{C}^{*})$ to the map $\overline{\beta}^{\vee}$.

Since $N$ is a finitely generated abelian group, we write it as
the invariant factor form
$$N=\mathbb{Z}^{d}\oplus\mathbb{Z}_{m_{1}}\oplus\cdots\oplus\mathbb{Z}_{m_{r}},$$
where $m_{1}|m_{2}|\cdots|m_{r}$. Let $\mathbf{\Sigma}$ be a
stacky fan. Suppose that the map $\beta$ generates the torsion
part of $N$, then $r\leq n-d$. We prove that the toric Deligne-Mumford stack $\mathcal{X}(\mathbf{\Sigma})$ is a nontrivial $\mu=\mu_{m_{1}}\times\cdots\times\mu_{m_{r}}$-gerbe
over the toric orbifold $\mathcal{X}(\mathbf{\Sigma_{red}})$
obtained as the stack of roots  of line bundles $M_{i}$. In
the Picard group of the toric orbifold, there exist $n-d$ line
bundles $M_{1},\cdots,M_{n-d}$ such that the $\mathcal{X}(\mathbf{\Sigma})$ can be constructed as $\mu_{m_{i}}$-gerbes over the line bundle $M_{i}$ for $i=1,\cdots,n-d$. These $n-d$ line bundles form the canonical generators of the Picard group.

Iwanari \cite{Iwa1}, \cite{Iwa2} proved that the integral Chow ring of the toric orbifold $\mathcal{X}(\mathbf{\Sigma_{red}})$ is isomorphic to the Stanley-Reisner ring $SR(\mathbf{\Sigma_{red}})$ of the reduced stacky fan $\mathbf{\Sigma_{red}}$, where
\begin{equation}\label{reduced-sr-ring}
SR(\mathbf{\Sigma_{red}}):=\frac{\mathbb{Z}[x_{i}:
\rho_{i}\in\Sigma(1)]}
{(I_{\Sigma}+C(\mathbf{\Sigma_{red}}))},
\end{equation}
the ideal $I_{\Sigma}$ is  generated by
\begin{equation}\label{ideal-i}
\{x_{i_{1}}\cdots x_{i_{k}}: \rho_{i_{1}}+\cdots+
\rho_{i_{k}}\notin\Sigma\},
\end{equation}
and $C(\mathbf{\Sigma_{red}})$ is the ideal generated by linear
relations:
\begin{equation}\Label{ideal-l-1}
\left(\sum_{i=1}^{n}\theta(b_{i})x_{i}\right)_{\theta\in N^{\star}}.
\end{equation}
In the toric orbifold case we have that
$$\frac{\mathbb{Z}[x_{i}:
\rho_{i}\in\Sigma(1)]}
{(C(\mathbf{\Sigma_{red}}))}\cong\overline{N}^{\vee},$$ where
$\overline{N}^{\vee}$ is the Picard group of the toric orbifold. Let
$t_{1},\cdots,t_{n-d}$ be the canonical generators of the Picard
group of $\mathcal{X}(\mathbf{\Sigma_{red}})$. Then $t_i$'s generate $x_i$'s
and $x_i$'s generate $t_i$'s, let $\mathbf{x}=(x_i),
\mathbf{t}=(t_i)$ be  column vectors, then there exist integral
matrices $A=(a_{i,j})$ and $C=(c_{i,j})$ such that
$\mathbf{x}=A\mathbf{t}$, and $\mathbf{t}=C\mathbf{x}$. Let
$\widetilde{\mathbf{x}}=(\widetilde{x}_i)$ be a column vector and 
$$
M=\left[\begin{array}{ccc}
m_1&&\\
&\ddots&\\
&&m_{n-d}\end{array}\right]
$$
a diagonal matrix such that 
\begin{equation}\label{matrix-key}
\widetilde{\mathbf{x}}:=AM\mathbf{t}=AMC\mathbf{x}=E\mathbf{x},
\end{equation}
where $E=AMC$ is a $n\times n$ integral matrix. So every
$\widetilde{x}_{i}$ is a integral linear combination of $x_{i}$'s.
Let $\mathcal{I}_{\mathbf{\Sigma}}$ be an ideal in
$\mathbb{Z}[x_{i}:\rho_{i}\in \Sigma(1)]$ generated by
\begin{equation}\label{ideal-product}
\{\widetilde{x}_{i_{1}}\cdots \widetilde{x}_{i_{k}}:
\rho_{i_{1}}+\cdots+ \rho_{i_{k}}\notin\Sigma\}.
\end{equation}
We define the Stanley-Reisner ring $SR(\mathbf{\Sigma})$ of the
stacky fan $\mathbf{\Sigma}$ as follows:
\begin{equation}\label{sr-ring}
SR(\mathbf{\Sigma}):=\frac{\mathbb{Z}[x_{i}: \rho_{i}\in\Sigma(1)]}
{(\mathcal{I}_{\mathbf{\Sigma}}+C(\mathbf{\Sigma}))},
\end{equation}
where $C(\mathbf{\Sigma})$ is the same as the ideal
$C(\mathbf{\Sigma_{red}})$ in (\ref{ideal-l-1}).

Let $A^{*}\left(\mathcal{X}(\mathbf{\Sigma}),\mathbb{Z}\right)$ be the integral
Chow ring of the toric Deligne-Mumford stack $\mathcal{X}(\mathbf{\Sigma})$. Then we have:
\begin{thm}\label{chowring-main1}
Let $\mathcal{X}(\mathbf{\Sigma})$ be a toric Deligne-Mumford stack with
semi-projective coarse moduli space. Suppose that in the map
$\beta$ of the stacky fan $\mathbf{\Sigma}$, the vectors $b_{1},\cdots,b_{n}$ generate
the torsion part of $N$. Then there is an isomorphism of rings:
$$A^{*}\left(\mathcal{X}(\mathbf{\Sigma}),\mathbb{Z}\right)\cong SR(\mathbf{\Sigma}).$$
\end{thm}
Since the toric Deligne-Mumford stack is a nontrivial $\mu$-gerbe constructed from a sequence of root gerbes of line bundles, Theorem \ref{chowring-main1} is proved by calculating Chow rings of root gerbes in terms of those of the bases and applying the result of Iwanari \cite{Iwa2}, see Section \ref{integralchow}.

The rational orbifold Chow ring of projective toric Deligne-Mumford stacks are computed in \cite{BCS}. Their result was generalized to semi-projective case in \cite{JT}. The orbifold Chow ring is isomorphic to the deformed ring of the fan. Using Theorem \ref{chowring-main1} we
compute the integral orbifold Chow ring for any semi-projective
toric Deligne-Mumford stacks.  For a stacky fan $\mathbf{\Sigma}$,
let $\sigma$ be a cone in $\Sigma$, define $Box(\sigma)$ to be all
the elements $v\in N$ such that $\overline{v}=
\sum_{\rho_{j}\subset\sigma}\alpha_{j}\overline{b}_{j}\in N$ for
$0\leq \alpha_{j}< 1$.  Let $Box(\mathbf{\Sigma})$ be the disjoint
union of all $Box(\sigma)$ for $\sigma\subset \Sigma$. The set $Box(\mathbf{\Sigma})$ is finite and may be listed as $\{v_1,\cdots,v_k\}$.

Now we consider the integral orbifold Chow ring of the toric
Deligne-Mumford stack $\mathcal{X}(\mathbf{\Sigma})$. First we introduce some notations.
Let $\widetilde{\mathbf{y}}=(\widetilde{y}^{b_{i}})$ and $\mathbf{y}=(y^{b_{i}})$ for $1\leq i\leq n$ be column vectors
such that 
$$\widetilde{\mathbf{y}}=E\mathbf{y},$$
where  the matrix $E=AMC$ is the same as the matrix in (\ref{matrix-key}). We
introduce the ring:
$$
S_{\mathbf{\Sigma}}:=\frac{\mathbb{Z}[y^{b_{i}}:
\rho_{i}\in\Sigma(1)]} {\mathcal{I}_{\mathbf{\Sigma}}},
$$
where $y$ is a formal variable and the ideal
$\mathcal{I}_{\mathbf{\Sigma}}$ is the ideal generated by the elements in (\ref{ideal-product}) replacing $\widetilde{x}_{i}$ by $\widetilde{y}^{b_{i}}$, i.e. 
by the elements
\begin{equation}\label{ideal-product2}
\{\widetilde{y}^{b_{i_{1}}}\cdots \widetilde{y}^{b_{i_{k}}}:
\rho_{i_{1}}+\cdots+ \rho_{i_{k}}\notin\Sigma\}.
\end{equation}
We define a ring
\begin{equation}\label{orbifoldchowring}
\mathbb{Z}[\mathbf{\Sigma}]=S_{\mathbf{\Sigma}}[y^{v_{1}},\cdots,y^{v_{k}}],
\end{equation}
which is a ring over the Stanley-Reisner ring. The product is
defined as follows: For any two $v_{1},v_{2}\in
Box(\mathbf{\Sigma})$, let
$v_{1}+v_{2}=\sum_{\rho_{j}\subset\sigma(\overline{v}_{1},\overline{v}_{2})}a_{j}b_{j}$
and let $I$ be the set of rays $\rho_{i}$ such that $a_{i}>1$, $J$
the set of rays $\rho_j$ such that
$\rho_j$ belongs to $\sigma(\overline{v}_{1}),\sigma(\overline{v}_{2})$, but
not $\sigma(\overline{v}_{3})$. Then
\begin{equation}\Label{product}
y^{v_{1}}\cdot y^{v_{2}}:=
\begin{cases}y^{\check{v}_{3}}\prod_{i\in I}\widetilde{y}^{b_{i}}\cdot
\prod_{i\in J}\widetilde{y}^{b_{i}}&\text{if
there is a cone}~ \sigma\in\Sigma ~\text{such that}~ \overline{v}_{1}\in\sigma, \overline{v}_{2}\in\sigma\,,\\
0&\text{otherwise}\,.\end{cases}
\end{equation}

Let $Cir(\mathbf{\Sigma})$ be the ideal  generated by the elements in (\ref{ideal-l-1}) replacing $x_{i}$ by $y^{b_{i}}$, i.e. 
by the linear relations:
\begin{equation}\Label{ideal-l}
\left(\sum_{i=1}^{n}\theta(b_{i})y^{b_{i}}\right)_{\theta\in
N^{\star}}.
\end{equation} 
Let
$A_{orb}^{*}\left(\mathcal{X}(\mathbf{\Sigma}),\mathbb{Z}\right)$ be
the integral orbifold Chow ring of the toric Deligne-Mumford stack
$\mathcal{X}(\mathbf{\Sigma})$. Then we have:
\begin{thm}\label{chowring-main2}
Let $\mathcal{X}(\mathbf{\Sigma})$ be a
toric Deligne-Mumford stack
associated to the
stacky fan $\mathbf{\Sigma}$ such that the coarse moduli space is semi-projective and the map $\beta$ generates the torsion part of $N$ in the stacky fan
$\mathbf{\Sigma}$. Then we have an isomorphism of graded rings:
$$A_{orb}^{*}\left(\mathcal{X}(\mathbf{\Sigma}),\mathbb{Z}\right)\cong \frac{\mathbb{Z}[\mathbf{\Sigma}]}{Cir(\mathbf{\Sigma})}.$$
\end{thm}
This is the first nontrivial examples of the formula for integral orbifold Chow rings.
Using Theorem \ref{chowring-main1}, the proof of Theorem \ref{chowring-main2}
is similar to \cite{BCS}, except that we work with integer coefficients.

This paper is organized as follows. In Section \ref{gerbe} we
recall the construction of Gale duality for finitely generated
abelian groups in \cite{BCS} and use it to study toric
Deligne-Mumford stacks. We give a construction of toric
Deligne-Mumford stacks in  this section. In Section
\ref{linebundle} we discuss line bundles over toric
Deligne-Mumford stacks and determine the Picard group of toric
Deligne-Mumford stacks.  In Section \ref{integralchow} we study
the integral Chow ring of toric Deligne-Mumford stacks. We prove
that the integral Chow ring of a toric Deligne-Mumford stack is
isomorphic to the Stanley-Reisner ring of its stacky fan. We study
the integral orbifold Chow ring in Section
\ref{orbifold-chowring}, and in Section \ref{example} we compute
some examples.
\subsection*{Conventions}
In this paper we work entirely algebraically over the field of
complex numbers. Chow rings and orbifold Chow rings are taken with
integer coefficients. By an orbifold we mean a smooth
Deligne-Mumford stack with trivial generic stabilizer.

We use $N^{\star}$ to denote the dual
$\text{Hom}_{\mathbb{Z}}(N,\mathbb{Z})$  and $\mathbb{C}^{*}$ the
multiplication group $\mathbb{C}-\{0\}$. We denote by $N\to
\overline{N}$ the natural map modulo torsion. Since $N$ is a
finitely generated abelian group, we write
$$N=\mathbb{Z}^{d}\oplus\mathbb{Z}_{m_{1}}\oplus\cdots\oplus\mathbb{Z}_{m_{r}},$$
where $m_{1}|m_{2}|\cdots|m_{r}$. This is called the invariant
factor decomposition.

\subsection*{Acknowledgments}
We thank Kai Behrend, Patrick Brosnan, Barbara Fantechi,  and Isamu Iwanari for valuable discussions.
Y. J. thanks the Fields Institute for Mathematical Science (Toronto, Canada) for financial
support during his visit in May--June, 2007. H.-H. T. thanks the Institut Mittag-Leffler (Djursholm, Sweden) for hospitality and support during his visit to the program ``moduli spaces''.
\section{Trivial and Nontrivial Gerbes}\label{gerbe}
In this section we first recall the construction of Gale duality for
finitely generated abelian groups. We use the properties of Gale
duality to classify toric Deligne-Mumford stacks as trivial and
nontrivial gerbes over the toric orbifolds.
\subsection{Gale duality and toric Deligne-Mumford
stacks.}\label{gale-toric} We recall the construction of Gale
duality according to \cite{BCS}. Let $N$ be a finitely generated
abelian group with rank $d$. Let
$$\beta: \mathbb{Z}^{n}\to N$$
be a map determined by $n$ integral vectors $\{b_1,\cdots,b_n\}$ in
$N$. Taking $\mathbb{Z}^{n}$ and $N$ as $\mathbb{Z}$-modules, from
the homological algebra, there exist projective resolutions
$\dot{E}$ and $\dot{F}$ of $\mathbb{Z}^{n}$ and $N$ satisfying the
following diagram
$$
\xymatrix{
~\dot{E}\rto\dto&~\mathbb{Z}^{n}\dto^{\beta}\\
~\dot{F}\rto&N.}
$$
Let $Cone(\beta)$ be the mapping cone of the map between $\dot{E}$
and $\dot{F}$. Then we have an exact sequence of the mapping cone:
$$0\longrightarrow\dot{F}\longrightarrow Cone(\beta)\longrightarrow\dot{E}[1]\longrightarrow 0,$$
where $\dot{E}[1]$ is the shifting of $\dot{E}$ by 1. Since
$\dot{E}$ is projective as $\mathbb{Z}$-modules, so we have the
exact sequence obtained by applying $\text{Hom}_\mathbb{Z}(-,\mathbb{Z})$:
$$0\longrightarrow\dot{E}[1]^{\star}\longrightarrow Cone(\beta)^{\star}\longrightarrow\dot{F}^{\star}\longrightarrow 0.$$
Taking cohomology of the above sequence we get the exact sequence
\begin{equation}\label{cohomology2}
N^{\star}\stackrel{\beta^{\star}}{\longrightarrow}(\mathbb{Z}^{n})^{\star}\stackrel{\beta^{\vee}}{\longrightarrow}
H^{1}(Cone(\beta)^{\star})\longrightarrow
Ext^{1}_{\mathbb{Z}}(N,\mathbb{Z})\longrightarrow 0.
\end{equation}

\begin{defn}\label{galedual}
Let $N^{\vee}=H^{1}(Cone(\beta)^{\star})$. The map
$$\beta^{\vee}: (\mathbb{Z}^{n})^{\star}\to N^{\vee}$$
is called the {\em Gale dual} of the map $\beta$. 
\end{defn} 
According to \cite{BCS}, both $N^{\vee}$ and $\beta^{\vee}$ are well defined up to natural isomorphism.

This construction can be made more clear. Since $N$ has rank $d$ and $\mathbb{Z}^{n}$ is a free $\mathbb{Z}$-module, the projection resolutions can be chosen as:
$$0\longrightarrow \mathbb{Z}^{n}\longrightarrow 0~=\dot{E},$$
$$0\longrightarrow\mathbb{Z}^{r}\stackrel{Q}{\longrightarrow}\mathbb{Z}^{d+r}\longrightarrow 0~=\dot{F},$$
where $Q$ is an integer matrix. Then there is a map $\mathbb{Z}^{n} \to \mathbb{Z}^{d+r}$ defined by a matrix $B$ which gives the map between $\dot{E}$ and $\dot{F}$. The mapping cone
$Cone(\beta)$ is given by the following complex:
$$0\longrightarrow\mathbb{Z}^{n+r}\stackrel{[B,Q]}{\longrightarrow}\mathbb{Z}^{d+r}\longrightarrow 0~=Cone(\beta).$$
(\ref{cohomology2}) is then obtained by applying the snake lemma to the following diagram
\begin{equation}\label{galedefinition}
\begin{CD}
0 @ >>>0@ >>> (\mathbb{Z}^{d+r})^{\star}@ >>>
(\mathbb{Z}^{d+r})^{\star} @
>>> 0\\
&& @VV{}V@VV{[B,Q]^{\star}}V@VV{}V \\
0@ >>> (Z^{n})^{\star} @ >{}>>(\mathbb{Z}^{n+r})^{\star}@ >>>
(\mathbb{Z}^{r})^{\star} @>>> 0.
\end{CD}
\end{equation}
Then $N^{\vee}=(\mathbb{Z}^{n+r})^{\star}/Im([B,Q]^{\star})$ and
$\beta^{\vee}$ is the composite map of the inclusion
$(\mathbb{Z}^{n})^{\star}\hookrightarrow(\mathbb{Z}^{n+r})^{\star}$
and the quotient map
$(\mathbb{Z}^{n+r})^{\star}\to(\mathbb{Z}^{n+r})^{\star}/Im([B,Q]^{\star})$.

\begin{rmk}\label{reducedcase}
If $N$ is free, i.e. there is no torsion part in the group $N$. Then
by (\ref{galedefinition}), the Gale dual $\beta^{\vee}$ is the
quotient map
$(\mathbb{Z}^{n})^{\star}\to (\mathbb{Z}^{n})^{\star}/Im([B]^{\star})$
and we have an exact sequence
$$N^{\star}\stackrel{\beta^{\star}}{\longrightarrow}(\mathbb{Z}^{n})^{\star}\stackrel{\beta^{\vee}}{\longrightarrow}
H^{1}(Cone(\beta)^{\star})\longrightarrow 0.$$
\end{rmk}

Let $N\to \overline{N}$ be the natural map of modding out torsion. Then $\overline{N}$ is a lattice. Let $\Sigma$ be a simplicial fan in the lattice $\overline{N}$ with $n$ rays $\{\rho_1,\cdots,\rho_n\}$. Choose $n$ integer vectors $\{b_1,\cdots,b_n\}$ such that $\overline{b}_{i}$ generates the ray $\rho_i$ for $1\leq i\leq n$. Then we have a map $\beta: \mathbb{Z}^{n}\to N$ determined by the vectors $\{b_1,\cdots,b_n\}$. We require that $\beta$ has finite cokernel.

\begin{defn}[\cite{BCS}]\label{stackyfan}
The triple $\mathbf{\Sigma}:=(N,\Sigma,\beta)$ is called a stacky fan.
\end{defn}

We define toric Deligne-Mumford stack from a stacky fan $\mathbf{\Sigma}$. Since $\beta$ has finite cokernel, by Proposition 2.2 and 2.3 in \cite{BCS}, we have the following exact sequences:
\begin{equation}\Label{exactchapter2-1}
0\longrightarrow
(N^{\vee})^{\star}\stackrel{(\beta^{\vee})^{\star}}{\longrightarrow}
\mathbb{Z}^{n}\stackrel{\beta}{\longrightarrow} N\longrightarrow
Coker(\beta)\longrightarrow 0,
\end{equation}
\begin{equation}\Label{exactchapter2-2}
0\longrightarrow N^{\star}\longrightarrow
\mathbb{Z}^{n}\stackrel{\beta^{\vee}}{\longrightarrow}
N^{\vee}\longrightarrow Coker(\beta^{\vee})\longrightarrow 0.
\end{equation}
Since $\mathbb{C}^{*}$ is a divisible as a $\mathbb{Z}$-module, applying $\text{Hom}_{\mathbb{Z}}(-,\mathbb{C}^{*})$ to (\ref{exactchapter2-2}) gives:

\begin{equation}\Label{exactchapter2-3}
1\longrightarrow \mu\longrightarrow
G\stackrel{\alpha}{\longrightarrow}
(\mathbb{C}^{*})^{n}\longrightarrow T\longrightarrow 1,
\end{equation}
where $\mu=\text{Hom}_{\mathbb{Z}}(Coker(\beta^{\vee}),\mathbb{C}^{*})$ is
finite, $G=\text{Hom}_{\mathbb{Z}}(N^{\vee},\mathbb{Z}^{*})$ and $T$ is the
$d$ dimensional torus $(\mathbb{C}^{*})^{d}$.

Let $\mathbb{C}[z_1,\cdots,z_n]$ be the coordinate ring of the
affine variety $\mathbb{A}^{n}$. Associated to the simplicial fan
$\Sigma$, there is an irrelevant ideal $J_{\Sigma}$ generated by the
elements:
\begin{equation}\label{idealchapter2}
\left\langle \prod_{\rho_{i}\nsubseteq \sigma}z_{i}: \sigma\in
\Sigma\right\rangle.
\end{equation}
Let $Z:=\mathbb{A}^{n}\setminus \mathbb{V}(J_{\Sigma})$. Then $Z$ is
a quasi-affine variety. The torus $(\mathbb{C}^{*})^{n}$ acts on $Z$
naturally since $Z$ is the complement of coordinate subspaces. The
algebraic  group $G$ acts on the variety $Z$ through the map
$\alpha$ in (\ref{exactchapter2-3}). Then we have
an action groupoid $Z\times G\toto Z$.

\begin{defn}[\cite{BCS}]\label{toricstackdefinition}
The toric Deligne-Mumford stack $\mathcal{X}(\mathbf{\Sigma})$
associated to the stacky fan $\mathbf{\Sigma}$ is defined to be the
quotient stack $[Z/G]$.
\end{defn}

Since $N$ is a finitely generated abelian group of rank $d$, we may write
$$N=\mathbb{Z}^{d}\oplus\mathbb{Z}_{m_{1}}\oplus\cdots\oplus\mathbb{Z}_{m_{r}}.$$
Then $\overline{N}=\mathbb{Z}^{d}$, and
$$\overline{\beta}:
\mathbb{Z}^{n}\to \overline{N}$$ is given by
$\{\overline{b}_{1},\cdots,\overline{b}_{n}\}$.  So
$\mathbf{\Sigma_{red}}:=(\overline{N},\Sigma,\overline{\beta})$ is a stacky fan. In the exact
sequence (\ref{exactchapter2-3}), let $\overline{G}=Im(\alpha)$,
then we have an exact sequence of abelian groups
$$1\longrightarrow \mu\longrightarrow G\longrightarrow \overline{G}\longrightarrow 1.$$
This is a central extension. By \cite{DP}, the quotient stack $[Z/G]$ is the
$\mu$-gerbe over the quotient stack $[Z/\overline{G}]=:\mathcal{X}(\mathbf{\Sigma_{red}})$ determined by this central extension. 

\begin{rmk}
The stack $\mathcal{X}(\mathbf{\Sigma_{red}})$ can be
constructed as follows. Consider  the following exact sequences
$$0\longrightarrow (\overline{N}^{\vee})^{\star}\longrightarrow
\mathbb{Z}^{n}\stackrel{\overline{\beta}}{\longrightarrow}
\overline{N}\longrightarrow 0;$$
where $\overline{\beta}$ is given by the vectors
$\{\overline{b}_1,\cdots,\overline{b}_n\}$, and
\begin{equation}\label{orbifold-exact}
0\longrightarrow N^{\star}\longrightarrow(\mathbb{Z}^{n})^{\star}\stackrel{\overline{\beta}^{\vee}}
\longrightarrow \overline{N}^{\vee}\longrightarrow 0.
\end{equation}
So
$A_{d-1}(X(\Sigma))=\overline{N}^{\vee}$ and from the construction
of Cox \cite{Cox},
$\mathcal{X}(\mathbf{\Sigma_{red}})=[Z/\overline{G}]$,
where $\overline{G}=\text{Hom}_{\mathbb{Z}}(\overline{N}^{\vee},\mathbb{C}^{*})$.
\end{rmk}

\subsection{Construction of toric Deligne-Mumford stacks.}
It is known that every toric Deligne-Mumford stack
$\mathcal{X}(\mathbf{\Sigma})$ is a $\mu$-gerbe over the
underlying toric orbifold for a finite abelian group $\mu$ and
some finite abelian gerbes over $\mathcal{X}(\mathbf{\Sigma})$ are
again toric Deligne-Mumford stacks, see \cite{Jiang2}. In this
section we classify trivial and nontrivial gerbes over toric
orbifolds.

\begin{lem}\label{smith-normal}
Let $\mathbb{Z}^{s}$ and $\mathbb{Z}^{t}$ be two free abelian
groups of ranks $s$ and $t$ respectively. Suppose that there is a
map $\beta: \mathbb{Z}^{s}\to\mathbb{Z}^{t}$ which
is given by an integral $t\times s$ matrix $A$. Then the dual map
$\beta^{\star}: (\mathbb{Z}^{t})^{\star}\to(\mathbb{Z}^{s})^{\star}$
is given by the transpose $A^{t}$ and
$$coker(\beta^{\star})\cong ker(\beta)\oplus coker{\beta}.$$
\end{lem}
\begin{pf}
For simplicity, we assume that $\beta$ has finite cokernel and
$s\geq t$. Since the matrix $A$ is an integer matrix, there exist
invertible $s\times s$ integer matrix $P$ and an $t\times t$
integer matrix $P^{'}$ such that $P^{'}AP$ is the matrix
$$\left[\begin{array}{cccc}
a_1&&\mathbf{0}&\mathbf{0} \\
~&\ddots & ~&\vdots\\
\mathbf{0}&&a_t&\mathbf{0}\\
\end{array}\right],
$$
with $a_1|\cdots|a_t$. This is the {\em Smith normal form} (see e.g. \cite{Pra}). It follows that 
$coker(\beta)\cong\mathbb{Z}_{a_{1}}\oplus\cdots\oplus\mathbb{Z}_{a_{t}}$.
After taking dual we get that the map
$\beta^\star: (\mathbb{Z}^{t})^{\star}\to(\mathbb{Z}^{s})^{\star}$
is given by the matrix
$$\left[\begin{array}{ccc}
a_1&&\mathbf{0} \\
~&\ddots & ~\\
\mathbf{0}&&a_t\\
\mathbf{0}&\cdots&\mathbf{0}
\end{array}\right].
$$
So it is easy to see that
$coker(\beta^{\star})=\mathbb{Z}^{s-t}\oplus\mathbb{Z}_{a_{1}}\oplus\cdots\oplus\mathbb{Z}_{a_{t}}$.
Since the kernel of $\beta$ is isomorphic to $\mathbb{Z}^{s-t}$,
we complete the proof.
\end{pf}

\begin{lem}\label{gen-tor}
Let $\mathbf{\Sigma}$ be a stacky fan. If the vectors
$\{b_{1},\cdots,b_{n}\}$ generate the torsion part
$\mathbb{Z}_{m_{1}}\oplus\cdots\oplus\mathbb{Z}_{m_{r}}$ of $N$,
then $r\leq n-d$ and $N^{\vee}\cong \overline{N}^{\vee}$.
\end{lem}
\begin{pf}
From the following diagram
$$\xymatrix{
~& ~\mathbb{Z}^{d+r}\dto \\
~\mathbb{Z}^{n}\urto\rto^{\beta}&~N,\dto\\
~&0 }
$$
It is easy to see that if $\beta$ generate the torision part of
$N$, then $r\leq n-d$.

By the construction of Gale duality in (\ref{galedefinition}),
$N^{\vee}\cong (\mathbb{Z}^{n+r})^{\star}/Im([B,Q]^{\star}).$
Since
$\overline{N}$ is free, we have 
$\overline{N}^{\vee}\cong
(\mathbb{Z}^{n})^{\star}/Im([\overline{B}]^{\star}).$ The map
$[B,Q]:
\mathbb{Z}^{n+r}\to\mathbb{Z}^{d+r}$ in the mapping
cone $Cone(\beta)$ has the same cokernel as the map $\beta$. Since
the map $\beta$ generate the torsion part of $N$, the map
$\beta$ has the same cokernel as the map $\overline{\beta}$. Thus
the map $[B,Q]$ has the same cokernel as the map
$[\overline{B}]:
\mathbb{Z}^{n}\to\mathbb{Z}^{d}.$ Also, $ker[B,Q]\simeq \mathbb{Z}^{n-d}\simeq ker[\overline{B}]$. By Lemma \ref{smith-normal}, we have $$N^\vee=coker [B,Q]^\star\simeq coker [\overline{B}]^\star=\overline{N}^\vee.$$
\end{pf}

Let $\mathbf{\Sigma}$ be the stacky fan and  $\mathbf{\Sigma_{red}}$
the corresponding reduced stacky fan.  Consider the following diagram
$$
\xymatrix{
~\mathbb{Z}^{n}\rto^{\beta}\dto_{id}&N\dto{}\\
~\mathbb{Z}^{n}\rto^{\overline{\beta}}&~\overline{N}.}
$$
Taking Gale dual yields
\begin{equation}\label{diagram-section4}
\begin{CD}
0 @ >>>\overline{N}^{\star}@ >>> \mathbb{Z}^{n}@ >{\overline{\beta}^{\vee}}>>
\overline{N}^{\vee} @
>>>0@>>> 0\\
&& @VV{}V@VV{}V@VV{\varphi}V@V{}VV \\
0@ >>> N^{\star} @ >{}>>\mathbb{Z}^{n}@ >{\beta^{\vee}}>> N^{\vee}@>>>cok(\beta^{\vee})@>>>0.
\end{CD}
\end{equation}

\begin{lem}\label{diagonal}
Let $\mathbf{\Sigma}=(N,\Sigma,\beta)$ be a stacky fan. If  the
vectors $\{b_{1},\cdots,b_{n}\}$ in the map $\beta$ generate the
torsion part
$\mathbb{Z}_{m_{1}}\oplus\cdots\oplus\mathbb{Z}_{m_{r}}$ of $N$,
then the map $\varphi$ in (\ref{diagram-section4}) is diagonalizable
over integers.
\end{lem}
\begin{pf}
First in the reduced stacky fan $\mathbf{\Sigma_{red}}$, the map
$\overline{\beta}: \mathbb{Z}^{n}\rightarrow \mathbb{Z}^{d}$ is
given by $\{\overline{b}_{1},\cdots,\overline{b}_{n}\}$. Consider
the following diagram
\begin{equation}\label{diagonal1}
\begin{CD}
~&&0@ >>> \mathbb{Z}^{n}@
>{id}>> \mathbb{Z}^{n} @
>>>0\\
&& @VV{\beta^{'}}V@VV{\beta}V@VV{\overline{\beta}}V \\
0@>>>\mathbb{Z}_{m_{1}}\oplus\cdots\oplus\mathbb{Z}_{m_{r}} @
>{}>>N@
>{}>> \overline{N}@>>>0.
\end{CD}
\end{equation}
From the definition of Gale dual in Section \ref{gale-toric}, we
have the morphisms of mapping cones $Cone(\beta^{'})$,
$Cone(\beta)$ and $Cone(\overline{\beta})$:
\begin{equation}\label{diagonal-1}
\begin{CD}
~ &&0 && 0 && ~ \\
&& @VV{}V@VV{}V  \\
0 @ >>>\mathbb{Z}^{r}@ >{\overline{Q}}>> \mathbb{Z}^{r}@>{}>> 0\\
&& @VV{}V@VV{}V \\
0@ >>> \mathbb{Z}^{n+r} @ >{[B,Q]}>>\mathbb{Z}^{d+r}@
>{}>>0\\
&& @VV{}V@VV{}V \\
0@ >>> \mathbb{Z}^{n} @ >{\overline{B}}>>\mathbb{Z}^{d}@
>>>0,\\
&& @VV{}V@VV{}V  \\
~ &&0 && 0 && ~
\end{CD}
\end{equation}
where $\overline{Q}$ is the diagonal matrix in $Q$. Dualizing gives the following diagram:
\begin{equation}\label{diagonal-2}
\begin{CD}
~ &&0 && 0 && 0&& ~ \\
&& @VV{}V@VV{}V@VV{}V  \\
0 @ >>>(\mathbb{Z}^{d})^{\star}@ >{\overline{B}^{\star}}>> (\mathbb{Z}^{n})^{\star}@>{\overline{\beta}^{\vee}}>>
\overline{N}^{\vee}@>>> 0\\
&& @VV{}V@VV{i}V@VV{\varphi}V \\
0@ >>> (\mathbb{Z}^{d+r})^{\star} @
>{[B,Q]^{\star}}>>(\mathbb{Z}^{n+r})^{\star}@
>{\pi}>>N^{\vee}@>>> 0\\
&& @VV{}V@VV{}V@VV{}V \\
0@ >>> (\mathbb{Z}^{r})^{\star} @
>{\overline{Q}^{\star}}>>(\mathbb{Z}^{r})^{\star}@
>>>\mathbb{Z}_{m_{1}}\oplus\cdots\oplus\mathbb{Z}_{m_{r}}@>>>0,\\
&& @VV{}V@VV{}V@VV{}V  \\
~ &&0 && 0 && 0&& ~
\end{CD}
\end{equation}
where $\pi\circ i=\beta^{\vee}$.  Since the map $\varphi$ is induced from the map
$i$ in (\ref{diagonal-2}), it is given
by an integer matrix $A$. Then from the general fact in the
finitely generated group theory there exist integer matrices $P$,
$P^{'}$ such that $PAP^{'}$ is a diagonal matrix with entries
$n_{1},\cdots,n_{s},0\cdots,0$ which satisfy the condition
$n_{1}|\cdots|n_{s}$. This is again the Smith normal form.
From the diagram (\ref{diagonal-2}) the third column is exact and
the cokernel is
$\mathbb{Z}_{m_{1}}\oplus\cdots\oplus\mathbb{Z}_{m_{r}}$, so the
diagonal matrix given by $\varphi$ is of the form
$$\left[\begin{array}{cccccc}
1&~&&&& \\
~&\ddots &&&& ~\\
&&1&&&\\
~&~&&m_{1}&&~\\
~&~&~&&\ddots &~\\
~&~&~&~&& m_{r} \end{array}\right].
$$
\end{pf}

\begin{prop}\label{nontrivialgerbe}
Let $\mathbf{\Sigma}=(N,\Sigma,\beta)$ be a stacky fan and $\mu=\mu_{m_{1}}\times\cdots\times\mu_{m_{r}}$. If  the vectors $\{b_{1},\cdots,b_{n}\}$ in the map $\beta$ generate the torsion part
$\mathbb{Z}_{m_{1}}\oplus\cdots\oplus\mathbb{Z}_{m_{r}}$ of $N$, then $\mathcal{X}(\mathbf{\Sigma})$ is a nontrivial $\mu$-gerbe over the toric orbifold $\mathcal{X}(\mathbf{\Sigma_{red}})$.
\end{prop}

\begin{pf}
Applying $\text{Hom}_{\mathbb{Z}}(-,\mathbb{C}^{*})$ to the 
diagram (\ref{diagram-section4}) yields
\begin{equation}\label{diagram-section4-2}
\begin{CD}
1 @ >>>\mu@ >>> G@ >{\alpha}>>
(\mathbb{C}^{*})^{n} @
>>>T@>>> 1\\
&& @VV{}V@VV{\alpha(\varphi)}V@VV{}V@V{}VV \\
1@ >>> 1 @ >{}>>\overline{G}@ >{\overline{\alpha}}>> (\mathbb{C}^{*})^{n}
@>>>T@>>> 1,
\end{CD}
\end{equation}
where the map $\alpha(\varphi)$ is given by the diagonal matrix:
\begin{equation}\label{diagonal-matrix}
\left[\begin{array}{cccccc}
(\cdot)^{1}&&&&\\
&\ddots&&&&\\
&&(\cdot)^{1}&&&\\
&&&(\cdot)^{m_{1}}&&\\
&&&&\ddots&\\
&&&&&(\cdot)^{m_r}\end{array}\right],
\end{equation}
since by Lemma \ref{diagonal} the map  
$\varphi$ in (\ref{diagram-section4}) is diagonalizable.
By Lemma \ref{gen-tor}, $N^{\vee}\cong \overline{N}^{\vee}$ and
so $G\cong \overline{G}$. So we get an exact sequence
\begin{equation}\label{central-extension}
1\longrightarrow\mu\stackrel{\varphi}{\longrightarrow}
G\longrightarrow\overline{G}\longrightarrow 1
\end{equation} which
is a central extension, where
$\mu=\text{Hom}_\mathbb{Z}(cok(\beta^{\vee}),\mathbb{C}^{*})=\mu_{m_{1}}\times\cdots\times\mu_{m_{r}}$. 
 We can decompose the left side of the diagram  (\ref{diagram-section4-2}) according to (\ref{diagonal-matrix}) to get the following diagram for each $\mu_{m_{i}}$:
\begin{equation}\label{diagram-section4-2-1}
\begin{CD}
1 @ >>>\mu_{m_i}@ >>> \mathbb{C}^{*}@ >{\alpha_i}>>
(\mathbb{C}^{*})^{n} @
>>>T@>>> 1\\
&& @VV{}V@VV{(\cdot)^{m_i}}V@VV{}V@V{}VV \\
1@ >>> 1 @ >{}>>\mathbb{C}^{*}@ >{\overline{\alpha}_i}>> (\mathbb{C}^{*})^{n}
@>>>T@>>> 1.
\end{CD}
\end{equation}
The corresponding component $\mathbb{C}^{*}\subset G$
determines a line bundle $M_{i}$ over the toric orbifold $\mathcal{X}(\mathbf{\Sigma_{red}})$. 
The central
extension (\ref{central-extension}) is nontrivial, because
$\overline{G}\cong G$ and the map $\varphi$ in
(\ref{central-extension}) is nontrivial on each component. So
from the definition of gerbes, the quotient stack
$\mathcal{X}(\mathbf{\Sigma})=[Z/G]$ is a nontrivial $\mu$-gerbe
over the toric orbifold
$\mathcal{X}(\mathbf{\Sigma_{red}})=[Z/\overline{G}]$ coming from the direct sum of line bundles $\oplus_{i}M_{i}$.
\end{pf}

\begin{rmk}\label{root-construction}
Recall that the stack of roots of a line bundle can be constructed as follows. More details can be found in Appendix B of \cite{AGV2}, or \cite{cad1}. Let $L$ be a line bundle over a variety (or an Artin stack) $X$. Let $m$ be a positive integer and consider the Kummer exact sequence:
$$1\longrightarrow \mu_{m}\longrightarrow \mathbb{C}^{*}\stackrel{(\cdot)^{m}}{\longrightarrow}\mathbb{C}^{*}\longrightarrow 1.$$
Then we have the quotient stack $[\mathbb{C}^{*}/\mathbb{C}^{*}]\simeq B\mu_m$, where the action is given by
$\lambda\cdot x=\lambda^{m}x$. Let $L^{*}$ be the complement of the zero section in the total space of $L$. The following twist
$$\sqrt[m]{L}:=[L^{*}/\mathbb{C}^{*}]=[L^{*}\times _{\mathbb{C}^{*}}\mathbb{C}^{*}/\mathbb{C}^{*}],$$
is the $\mu_{m}$-gerbe over $X$ coming from the line bundle $L$, which is called the stack of $m$-th root of $L$. $\sqrt[m]{L}$ may be viewed as a toric stack bundle in the sense of
\cite{Jiang} since the stack $[\mathbb{C}^{*}/\mathbb{C}^{*}]$ is a toric Deligne-Mumford stack.

In the proof of Proposition \ref{nontrivialgerbe}, the toric Deligne-Mumford stack $\mathcal{X}(\mathbf{\Sigma})$ is obtained by applying a sequence of root constructions to line bundles $M_{i}$ over the toric orbifold  $\mathcal{X}(\mathbf{\Sigma_{red}})$, i.e. 
$$\mathcal{X}(\mathbf{\Sigma})\cong \sqrt[m_{1}]{M_{1}}\times_{\mathcal{X}(\mathbf{\Sigma_{red}})}\cdots\times_{\mathcal{X}(\mathbf{\Sigma_{red}})}  
\sqrt[m_{r}]{M_{r}}.$$  

Similar descriptions have also been obtained independently by F. Perroni \cite{Perroni}.
\end{rmk}

\begin{example}
Let $N=\mathbb{Z}^{2}\oplus\mathbb{Z}_{2}\oplus\mathbb{Z}_{4}$,
and $\beta: \mathbb{Z}^{4}\to N$ be given by the vectors
$$\{b_{1}=(1,0,1,0),b_{2}=(0,1,0,0),b_{3}=(-1,2,0,0),b_{4}=(0,-1,0,1)\}.$$
Then $\overline{N}=\mathbb{Z}^{2}$ and let $\Sigma$ be the
complete fan in $\mathbb{R}^{2}$ generated by
$\{\overline{b}_{1},\overline{b}_{2},\overline{b}_{3},\overline{b}_{4}\}$.
The fan $\Sigma$ is  the fan of Hirzebruch surface
$\mathbb{F}_{2}$. Then $\mathbf{\Sigma}=(N,\Sigma,\beta)$ is a
stacky fan. We compute that $\beta^{\vee}:
\mathbb{Z}^{4}\to N^{\vee}=\mathbb{Z}^{2}$ is given by
the matrix
$$\left[
\begin{array}{cccc}
  2&-4&2&0 \\
  0&4&0&4\\
\end{array}
\right].$$ The reduced stacky fan is
$\mathbf{\Sigma_{red}}=(\mathbb{Z}^{2},\Sigma,\overline{\beta})$,
where $\overline{\beta}:\mathbb{Z}^{4}\to
\mathbb{Z}^{2}$ is given by
$\{\overline{b}_{1},\overline{b}_{2},\overline{b}_{3},\overline{b}_{4}\}$.
The Gale dual $\overline{\beta}^{\vee}:
\mathbb{Z}^{4}\to N^{\vee}=\mathbb{Z}^{2}$ is given by
the matrix
$$\left[
\begin{array}{cccc}
  1&-2&1&0 \\
  0&1&0&1\\
\end{array}
\right].$$ So in this example the diagram (\ref{diagram-section4})
is:
\begin{equation}\label{diagram-section4-example}
\begin{CD}
0 @ >>>\mathbb{Z}^{2}@ >>> \mathbb{Z}^{4}@
>{\overline{\beta}^{\vee}}>> \mathbb{Z}^{2} @
>>>0@>>> 0\\
&& @VV{}V@VV{id}V@VV{\varphi}V@V{}VV \\
0@ >>> \mathbb{Z}^{2} @ >{}>>\mathbb{Z}^{4}@
>{\beta^{\vee}}>>
\mathbb{Z}^{2}@>>>\mathbb{Z}_{2}\oplus\mathbb{Z}_{4}@>>>0,
\end{CD}
\end{equation}
where $\varphi$ is the diagonal matrix $$\left[
\begin{array}{cc}
  2&0 \\
  0&4\\
\end{array}
\right].$$ So from (\ref{diagram-section4-example}) we get the
following exact sequence:
\begin{equation}\label{gerbe-example}
\begin{CD}
1 @ >>>\mu_{2}\times \mu_{4}@ >>> (\mathbb{C}^{*})^{2}@
>{}>> (\mathbb{C}^{*})^{4} @
>>>(\mathbb{C}^{*})^{2}@>>> 1\\
&& @VV{}V@VV{\alpha(\varphi)}V@VV{id}V@V{}VV \\
1@ >>> 1 @ >{}>>(\mathbb{C}^{*})^{2}@
>{}>>
(\mathbb{C}^{*})^{4}@>>>(\mathbb{C}^{*})^{2}@>>>1.
\end{CD}
\end{equation}
The toric Deligne-Mumford stack
$\mathcal{X}(\mathbf{\Sigma})=[\mathbb{C}^{4}-\mathbb{V}(J_{\Sigma})/(\mathbb{C}^{*})^{2}]$,
where the action is given by the transpose of the matrix
$\beta^{\vee}$. Let $L_{1},L_{2}$ be the two line bundles over
$\mathbb{F}_{2}$ which are the two canonical generators in the
Picard group $\mathbb{Z}^{2}$ of $\mathbb{F}_{2}$. The toric
Deligne-Mumford stack $\mathcal{X}(\mathbf{\Sigma})$ is a
$\mu_{2}\times\mu_{4}$-gerbe over the Hirzebruch surface
$\mathcal{X}(\mathbf{\Sigma_{red}})=\mathbb{F}_{2}$ coming from
the direct sum of line bundles $L_{1}\oplus L_{2}$. $\mathcal{X}(\mathbf{\Sigma})$ can be
constructed by taking square root and quartic root of line bundles $L_{1}$
and $L_{2}$ respectively: $\mathcal{X}(\mathbf{\Sigma})\cong \sqrt[2]{L_{1}}\times_{\mathbb{F}_{2}}\sqrt[4]{L_{2}}.$
\end{example}

\begin{rmk}\label{rmk-zero}
In the map $\beta: \mathbb{Z}^{n}\to N$, if the
components of the torsion part of $b_{i}$ are zero, then it is
easy to check that
$N^{\vee}\cong\overline{N}^{\vee}\oplus\mathbb{Z}_{m_{1}}\oplus\cdots\oplus\mathbb{Z}_{m_{r}}$
and
$\mathcal{X}(\mathbf{\Sigma})=\mathcal{X}(\mathbf{\Sigma_{red}})\times\mathcal{B}\mu$
with $\mu=\mu_{m_{1}}\times\cdots\times\mu_{m_{r}}$.
\end{rmk}

\begin{prop}\label{toric-gerbe}
Every toric Deligne-Mumford stack $\mathcal{X}(\mathbf{\Sigma})$ has
a  decomposition:
$$\mathcal{X}(\mathbf{\Sigma})\cong \mathcal{X}(\mathbf{\Sigma^{'}})\times \mathcal{B}\mu,$$
where $\mathcal{X}(\mathbf{\Sigma^{'}})$ is a nontrivial gerbe
over the toric orbifold $\mathcal{X}(\mathbf{\Sigma_{red}^{'}})$
and $\mathcal{B}\mu$ is the classifying stack for a finite abelian
group $\mu$.
\end{prop}
\begin{pf}
Consider the map $\beta: \mathbb{Z}^{n}\to N$ given by integral vectors $\{b_{1},\cdots,b_{n}\}\subset N$.
Let $N=\mathbb{Z}^{d}\oplus\mathbb{Z}_{m_{1}}\oplus\cdots\oplus\mathbb{Z}_{m_{r}}$
and $N_{tor}=\mathbb{Z}_{m_{1}}\oplus\cdots\oplus\mathbb{Z}_{m_{r}}$.
Then there exists a subgroup $N_{tor}^{'}\subset N_{tor}$ such that that $\{b_{1},\cdots,b_{n}\}$
generate  $N_{tor}^{'}$. Let $N^{'}=\mathbb{Z}^{d}\oplus N^{'}_{tor}$ and let
$\beta^{'}: \mathbb{Z}^{n}\to N^{'}$ be the map given by
$\{b_{1},\cdots,b_{n}\}$.
Then $\mathbf{\Sigma^{'}}=(N^{'},\Sigma,\beta^{'})$ is a stacky fan such that
$\{b_{1}^{'},\cdots,b_{n}^{'}\}$ generate the torsion part of $N^{'}$.

Let $\mathbf{\Sigma_{red}}$ be the corresponding reduced stacky fan. Then by Lemma \ref{gen-tor},
in the Gale dual $(\beta^{'})^{\vee}:\mathbb{Z}^{n}\to (N^{'})^{\vee}$ and
$\overline{\beta}^{\vee}:\mathbb{Z}^{n}\to (\overline{N})^{\vee}$, 
$(N^{'})^{\vee}\cong (\overline{N})^{\vee}$. Since the cokernel of the map
$\beta$ is 
$$cok(\overline{\beta})\oplus N_{tor}/N_{tor}^{'}\cong 
cok(\beta^{'})\oplus N_{tor}/N_{tor}^{'},$$ 
by Lemma \ref{smith-normal},
the Gale dual map of $\beta$ is
$$\beta^{\vee}: \mathbb{Z}^{n}\to (N^{'})^{\vee}\oplus N_{tor}/N^{'}_{tor}.$$
Let $\mu=N_{tor}/N^{'}_{tor}$.
By Lemma \ref{nontrivialgerbe}
and Remark \ref{rmk-zero} the proposition follows.
\end{pf}

\section{Line Bundles over Toric Deligne-Mumford Stacks}\label{linebundle}

In this section we prove that the Picard group of a toric Deligne-Mumford stacks is isomorphic to
$N^{\vee}$ in the Gale dual $\beta^{\vee}: \mathbb{Z}^{n}\rightarrow N^{\vee}$ of $\beta: \mathbb{Z}^{n}\rightarrow N$.

Let $\mathcal{X}(\mathbf{\Sigma})=[Z/G]$ be a toric Deligne-Mumford stack associated to the
stacky fan $\mathbf{\Sigma}=(N,\Sigma,\beta)$. In the case of a quotient stack, a line bundle $\mathcal{L}$ on $\mathcal{X}(\mathbf{\Sigma})$ is a $G$-equivariant line bundle $L$ on $Z$. A $G$-equivariant line bundle on $Z$ is a line bundle $L$ on $Z$ together with an isomorphism $\varphi: pr^{*}L\rightarrow\mu^{*}L$, where $pr: G\times Z\rightarrow Z$ is the projection and $\mu: G\times Z\rightarrow Z$ is the action. Let $Pic(\mathcal{X}(\mathbf{\Sigma}))$ denote the Picard group of $\mathcal{X}(\mathbf{\Sigma})$.

\begin{prop}\label{picardgroup}
$Pic(\mathcal{X}(\mathbf{\Sigma}))\cong N^{\vee}$.
\end{prop}
\begin{pf}
Since $G=\text{Hom}_\mathbb{Z}(N^\vee, \mathbb{C}^*)$, Pontryagin duality implies that the character group $G^\vee$ of $G$ is isomorphic to $N^\vee$.

As discussed above, $\mathcal{X}(\mathbf{\Sigma})=[Z/G]$ where $Z=\mathbb{A}^n\setminus \mathbb{V}(J_\Sigma)$ with $\text{codim}\,\mathbb{V}(J_\Sigma)\geq 2$. This implies that $Pic(Z)\simeq Pic(\mathbb{A}^n)\simeq \mathbb{Z}$, generated by the trivial line bundle.

Let $\mathcal{L}\in Pic (\sX(\mathbf{\Sigma})$, then $\sL$ corresponds to a $G$-equivariant line bundle $L\to Z$, which is trivial but not necessarily $G$-equivariantly trivial. The $G$-action on $L$ is given by a character $\chi_\sL: G\to \mathbb{C}^*$. Thus there is a map $$Pic(\sX(\mathbf{\Sigma}))\to G^\vee, \quad \sL\to \chi_\sL.$$
On the other hand, any character $\chi: G\to \mathbb{C}^*$ gives a $G$-equivariant line bundle $L_\chi$ over $Z$, hence a line bundle on $\sX(\mathbf{\Sigma})$. The result follows.
\end{pf}

From the results in \cite{EG}, \cite{Kr} that the Picard group of a quotient  stack
is isomorphic to its first Chow group, we have:
\begin{cor}\label{picard-chow}
There is an isomorphism $A^{1}(\mathcal{X}(\mathbf{\Sigma}),\mathbb{Z})\cong N^{\vee}$.
\end{cor}


For each ray $\rho_{i}$, we define a line bundle $\mathcal{L}_{i}$ over
$\mathcal{X}(\mathbf{\Sigma})$ to be the trivial line bundle $\mathbb{C}$ with the
$G$ action on $\mathbb{C}$ given by the $i$-th component of the map
$\alpha: G\to (\mathbb{C}^{*})^{n}$.
We can similarly define a line bundle $L_{i}$ over the
toric orbifold $\sX(\mathbf{\Sigma_{red}})$. This line bundle
$L_{i}$ is the trivial line bundle $\mathbb{C}$ over $Z$ on which the action of $\overline{G}$ is through the $i$-th component of the map
$\overline{\alpha}: \overline{G}\to
(\mathbb{C}^{*})^{n}$, where $\overline{\alpha}$ is obtained by
taking $\text{Hom}_\mathbb{Z}(-,\mathbb{C}^{*})$ to the map $\overline{\beta}^{\vee}$
in (\ref{orbifold-exact}). 

For the stacky fan $\mathbf{\Sigma}$  and its reduced stacky fan
$\mathbf{\Sigma_{red}}$, we consider again the diagram 
(\ref{diagram-section4}).
Since we have that
$cok(\beta^{\vee})=\mathbb{Z}_{m_{1}}\oplus\cdots\oplus\mathbb{Z}_{m_{n-d}}$ by Lemma \ref{diagonal}, the map $\varphi$ is given by the diagonal matrix
$$\left[\begin{array}{ccc}
m_{1}&~&~\\
~&\ddots&~\\
~&~&m_{n-d} \end{array}\right],
$$
where some $m_{i}$'s are $1$ since $r\leq n-d$.

Let $\mathbf{t}:=\{t_{1},\cdots,t_{n-d}\}$ be the generators of $\overline{N}^{\vee}$ so that the map $\varphi$ is diagonalized. Let $\{e_1,...,e_n\}$ be the standard basis of $\mathbb{Z}^n$. Let
$\mathbf{x}=(x_i)$ and $\mathbf{t}=(t_i)$ be column vectors, then there exist a matrix $A$
such that $\mathbf{x}=A\mathbf{t}$. Then under these bases the map $\overline{\beta}^{\vee}$ is given by the matrix $A^t$. Suppost that the map
$\beta^{\vee}$ is given by a matrix $B$, then we have $B=MA^t$. Since $\overline{\beta}^\vee$ is surjective, there exists an
integral matrix $C$ such that $\mathbf{t}=C\mathbf{x}$, where
$\mathbf{x}:=(x_1,\cdots,x_n)^t$ with $x_j:= Ae_j$. Let
$\widetilde{\mathbf{x}}:=(\widetilde{x}_1,\cdots,\widetilde{x}_n)^t$
be defined by 
\begin{equation}\label{key-formula}
\widetilde{\mathbf{x}}=AM\mathbf{t}=AMC\mathbf{x}.
\end{equation}
Then every $\widetilde{x}_i$ is an integral linear combination of
$x_{i}$'s by the above formula. We call the formula
(\ref{key-formula}) the associated formula of the stacky fan
$\mathbf{\Sigma}$.

Let $\pi: \mathcal{X}(\mathbf{\Sigma})\to
\mathcal{X}(\mathbf{\Sigma_{red}})$ be the rigidification map. For the line bundle $L_j$ over
$\mathcal{X}(\mathbf{\Sigma_{red}})$ corresponding to ray $\rho_j$, we have the pullback
$\pi^{*}(L_j)=\mathcal{L}_j$ over $\mathcal{X}(\mathbf{\Sigma})$.
From Lemma \ref{diagonal},  we have the morphism 
$\varphi: \overline{N}^{\vee}\longrightarrow N^{\vee}$ after choosing  the basis
of $\overline{N}^{\vee}\cong N^{\vee}$ which is diagonalizable. Let $id: \overline{N}^{\vee}\longrightarrow N^{\vee}$ be the 
identity morphism under the basis, dualizing we an isomorphism  
$$id: G\stackrel{\cong}{\longrightarrow} \overline{G}.$$
Under this isomorphism, we can  take $L_j$ as a line bundle over $\mathcal{X}(\mathbf{\Sigma})$ using the 
character
$\overline{\alpha}_j: \overline{G}\to \mathbb{C}$ in the 
$j$-th component of the map $\overline{\alpha}$.
Suppose that in the formula (\ref{key-formula}),
$$\widetilde{x}_i=\sum_{j=1}^{n}a_{j,i}x_{j},$$
where $a_{1,i},\cdots,a_{n,i}$ are integers, then we have the
following proposition.
\begin{prop}\label{linebundle-com}
$\mathcal{L}_{i}=\bigotimes_{j}L_{j}^{\otimes a_{j,i}}$.
\end{prop}
\begin{pf}
From the construction of the line bundle $\mathcal{L}_{i}$ on the
toric Deligne-Mumford stack the proposition can be easily proved
from the diagram (\ref{diagram-section4}) and the formula
(\ref{key-formula}).
\end{pf}
\section{Integral Chow Ring of Toric Deligne-Mumford Stacks}\label{integralchow}

In this section we study the integral Chow ring of toric Deligne-Mumford stacks. For references of integral Chow ring of stacks, see \cite{EG} and \cite{Kr}. In this section, every Deligne-Mumford stack $\mathcal{X}(\mathbf{\Sigma})$ is semi-projective satisfying the condition of Lemma \ref{nontrivialgerbe}. We use the results in \cite{Iwa1}, \cite{Iwa2} concerning the integral Chow ring of a simplicial toric orbifold.

\subsection{The Chow ring of stack of roots of a line bundle}

In this section we compute the Chow ring with integer coefficients of the stack $\sqrt[m]{L}$ of $m$-th root of a line bundle $L$.

Let $\mathbb{C}^*$ act on the total space of $L$ by acting trivially on the base $X$ and acting by $\lambda\cdot x:= \lambda^m x$ on the fiber. As recalled in Remark \ref{root-construction}, $\sqrt[m]{L}$ is the stack quotient $[L^*/\mathbb{C}^*]$ with respect to this $\mathbb{C}^*$-action. 

Let $i: X\hookrightarrow L$ be the inclusion of the zero section, and  $j: L\setminus X\hookrightarrow L$ the inclusion of its complement. Then we have an exact sequence on the $\mathbb{C}^{*}$-equivariant
Chow groups:
\begin{equation}\label{equivariant-exact}
A_{*}(X)_{\mathbb{C}^{*}}\stackrel{i_{*}}{\longrightarrow}A_{*}(L)_{\mathbb{C}^*}\stackrel{j^{*}}{\longrightarrow}
A_{*}(L\setminus X)_{\mathbb{C}^*}\longrightarrow 0,
\end{equation}
where $i_{*}$ is the pushforward and $j^{*}$ is the flat pullback.
By \cite{EG},
$$A^{*}([L^{*}/\mathbb{C}^{*}])=A_{\mathbb{C}^{*}}^{*}(L\setminus X).$$
Let $\pi: L\to X$ be the structure map of the line bundle. Then we have
$$\pi^{*}: A_{*}(X)_{\mathbb{C}^{*}}\stackrel{\cong}{\longrightarrow}A_{*+1}(L)_{\mathbb{C}^{*}}.$$
Let
\begin{equation}\label{mapb}
b:=(\pi^{*})^{-1}\circ i_*: A_{*}(X)_{\mathbb{C}^{*}}\to A_{*-1}(X)_{\mathbb{C}^{*}}.
\end{equation}
We may rewrite the sequence (\ref{equivariant-exact}) as
\begin{equation}\label{equivariant-exact2}
A_{*}(X)_{\mathbb{C}^{*}}\stackrel{b}{\longrightarrow}A_{*-1}(X)_{\mathbb{C}^{*}}\stackrel{j^{*}}{\longrightarrow}
A_{*}(L\setminus X)_{\mathbb{C}^{*}}\longrightarrow 0.
\end{equation}
This implies that 
\begin{equation}\label{Chow-root}
A_{*}([L^{*}/\mathbb{C}^{*}])\cong A_{*-1}(X)_{\mathbb{C}^{*}}/Im(b).
\end{equation}
Clearly the map $b$ is given by the $\mathbb{C}^{*}$-equivariant first Chern class of the line bundle
$L$, which is $c_{1}(L)-mt,$ where $c_{1}(L)$ is the non-equivariant first Chern class of $L$, and
$t$ is the equivariant parameter. Thus we have the following
proposition:
\begin{prop}\label{Chow-root2}
The Chow ring $A^{*}([L^{*}/\mathbb{C}^{*}])$ is isomorphic to the quotient ring
$$\frac{A^{*}(X)[t]}{(c_{1}(L)-mt)}.$$
\end{prop}
\begin{pf}
Since $\mathbb{C}^*$ acts trivially on $X$, we have $A^*(X)_{\mathbb{C}^*}=A^*(X)[t]$. Since $c_1(L)-mt$ is not a zero divisor, (\ref{equivariant-exact2}) is exact on the left, and the image of $b$ is the ideal generated by $c_1(L)-mt$.
\end{pf}

\subsection{Integral Chow ring of toric Deligne-Mumford stacks.}

Let $\mathbf{\Sigma}$ be a stacky fan such that
the map $\beta$ generate the torsion part of the group $N$ and  $\mathcal{X}(\mathbf{\Sigma})$
the associated toric Deligne-Mumford stack . Let
$\mathcal{X}_{orb}(\Sigma)$ be the underlying toric orbifold given
by the simplicial fan $\Sigma$.
Let
\begin{equation}\label{stanley-orbifold}
\frac{\mathbb{Z}[x_{i}: \rho_{i}\in\Sigma(1)]}{\left(I_{\Sigma}+Cir(\Sigma)\right)}
\end{equation}
be the Stanley-Reisner ring of the fan $\Sigma$, where $x_{i}$ corresponds to the  torus invariant
divisor $D_{\rho_{i}}$, $Cir(\Sigma)$ is the ideal generated by the
linear relations:
$$\left(\sum_{\rho_{i}\in\Sigma}\theta(v_{i})x_{i}\right)_{\theta\in N^{\star}},$$
where $v_{i}$ is the first lattice point of the ray $\rho_{i}$, and $I_{\Sigma}$ is the ideal generated by
square-free monomials
$$\{ x_{i_{1}}\cdots x_{i_{k}}: \rho_{i_{1}}+\cdots +\rho_{i_{k}}~\mbox{is not a cone}~\sigma\in
\Sigma\}$$
in (\ref{ideal-i}).

\begin{prop}[\cite{Iwa1}]\label{integral-orbifold}
The integral Chow ring $A^{*}(\mathcal{X}_{orb}(\Sigma),\mathbb{Z})$ is isomorphic to (\ref{stanley-orbifold}).
\end{prop}

Consider $\mathbf{\Sigma_{red}}=(\overline{N},\Sigma,\overline{\beta})$, where
$\overline{\beta}: \mathbb{Z}^{n}\to \overline{N}$ is given by the vectors
$\{\overline{b}_{1},\cdots,\overline{b}_{n}\}$, and the toric Deligne-Mumford stack
$\mathcal{X}(\mathbf{\Sigma_{red}})=[Z/\overline{G}]$. Let $C(\mathbf{\Sigma_{red}})$ is the ideal generated by the
linear relations:
$$\left(\sum_{\rho_{i}\in\Sigma}\theta(b_{i})x_{i}\right)_{\theta\in N^{\star}}$$
in (\ref{ideal-l-1}).

\begin{prop}[\cite{Iwa2}]\label{integral-orbifold2}
There is an isomorphism of rings:
$$A^{*}(\mathcal{X}(\mathbf{\Sigma_{red}}),\mathbb{Z})\cong
\frac{\mathbb{Z}[x_{i}:
\rho_{i}\in\Sigma(1)]}{\left(I_{\Sigma}+C(\mathbf{\Sigma_{red}})\right)}.$$
\end{prop}

\begin{rmk}
Let $\mathbf{\Sigma}=(N,\Sigma,\beta)$ be a stacky fan. For the
simplicial fan $\Sigma$, the toric orbifold
$\mathcal{X}_{orb}(\Sigma)$ associated to $\Sigma$ has stack structures in codimension at least $2$. The toric orbifold
$\mathcal{X}(\mathbf{\Sigma_{red}})$ has stack structures in codimension
at least $1$. The toric orbifold
$\mathcal{X}(\mathbf{\Sigma_{red}})$ can be obtained from
$\mathcal{X}_{orb}(\Sigma)$ by taking roots of divisors. Since we don't need this result here, we omit the
details.
\end{rmk}

In view of Proposition \ref{nontrivialgerbe}, the idea of proving Theorem \ref{chowring-main1} is to compute the integral Chow ring of $\mathcal{X}(\mathbf{\Sigma})$ by combining Propositions \ref{integral-orbifold2} and \ref{Chow-root2}.


\subsection*{Proof of Theorem \ref{chowring-main1}.}
As in the proof of Proposition \ref{nontrivialgerbe}, let $M_{i}\to \mathcal{X}(\mathbf{\Sigma_{red}})$ be the line bundle over the toric orbifold given by one generator in
$\overline{N}^{\vee}$ and Let $m_{i}$ be the corresponding positive
integer. The toric Deligne-Mumford stack $\mathcal{X}(\mathbf{\Sigma})$ is a nontrivial
$\mu_{m_{1}}\times\cdots\times\mu_{m_{r}}$-gerbe over the toric
orbifold $\mathcal{X}(\mathbf{\Sigma_{red}})$ obtained from a
sequence of root gerbe constructions determined by the line
bundles $M_{i}$ for $1\leq i\leq n-d$.

Since the integral Chow ring
$A^{*}(\mathcal{X}(\mathbf{\Sigma_{red}}),\mathbb{Z})$ is generated
by the Picard group $\overline{N}^{\vee}$ and
$\frac{\mathbb{Z}[x_{i}:
\rho_{i}\in\Sigma(1)]}{C(\mathbf{\Sigma_{red}})}\cong \overline{N}^{\vee}$, so by Proposition \ref{integral-orbifold2}, 
\begin{equation}
A^{*}(\mathcal{X}(\mathbf{\Sigma_{red}});\mathbb{Z}) \cong \frac{\mathbb{Z}[x_{i}:
\rho_{i}\in\Sigma(1)]}{(\mathcal{I}_{\mathbf{\Sigma_{red}}}+C(\mathbf{\Sigma_{red}}))}
\nonumber 
\cong
\frac{\mathbb{Z}[t_{1},\cdots,t_{n-d}]}{I_{\mathbf{\Sigma_{red}}}},
\nonumber
\end{equation}
where $\{t_1,...,t_{n-d}\}$ is a basis of $N^\vee$, the ideal $I_{\mathbf{\Sigma_{red}}}$ is obtained from
$\mathcal{I}_{\mathbf{\Sigma_{red}}}$ by expressing $x_{i}^{,}$s in terms of $t_{i}^{,}$s. By Proposition \ref{Chow-root2}, the Chow ring
$A^{*}(\mathcal{X}(\mathbf{\Sigma}),\mathbb{Z})$  is isomorphic to the ring obtained from
$\frac{\mathbb{Z}[t_{1},\cdots,t_{n-d}]}{I_{\mathbf{\Sigma_{red}}}}$
by replacing the canonical generators $\{t_1,\cdots,t_{n-d}\}$ by
$\{m_1t_1,\cdots,m_{n-d}t_{n-d}\}.$
In view of (\ref{key-formula}), this  is isomorphic to the
Stanley-Reisner ring $SR(\mathbf{\Sigma})$ of the stacky fan
$\mathbf{\Sigma}$ since in the ring $SR(\mathbf{\Sigma})$, the ideal
$\mathcal{I}_{\mathbf{\Sigma}}$ is obtained from the ideal
$\mathcal{I}_{\mathbf{\Sigma_{red}}}$ in (\ref{ideal-product}) replacing $x_i$ by
$\widetilde{x}_{i}$ for each ray $\rho_{i}$. $\square$

By Proposition \ref{toric-gerbe},
every toric Deligne-Mumford stack $\mathcal{X}(\mathbf{\Sigma})$ has
a  decomposition: $\mathcal{X}(\mathbf{\Sigma})\cong \mathcal{X}(\mathbf{\Sigma^{'}})\times \mathcal{B}\mu,$
where $\mathcal{X}(\mathbf{\Sigma^{'}})$ is a nontrivial gerbe
over the toric orbifold $\mathcal{X}(\mathbf{\Sigma_{red}})$ and
$\mathcal{B}\mu$ is the classifying stack of a finite abelian
group $\mu=\mu_{r_1}\times\cdots\times\mu_{r_{s}}$.  It is known that $$A^{*}(\mathcal{B}\mu,\mathbb{Z})\cong
\mathbb{Z}[t_1,\cdots,t_{s}]/(r_1t_1,\cdots,r_st_s).$$
Thus we have 
\begin{prop}\label{Chow-stack}
The integral Chow ring of $\mathcal{X}(\mathbf{\Sigma})$ is given by 
$$A^{*}(\mathcal{X}(\mathbf{\Sigma}),\mathbb{Z})\cong
A^{*}(\mathcal{X}(\mathbf{\Sigma^{'}}),\mathbb{Z})[t_1,\cdots,t_{s}]/(r_1t_1,\cdots,r_st_s).$$ 
\end{prop}
\section{The Integral Orbifold Chow Ring}\label{orbifold-chowring}

In this section we compute the integral orbifold Chow ring of
toric Deligne-Mumford stacks.

\subsection{The inertia stack.}

Let $\mathcal{X}(\mathbf{\Sigma})$ be a toric Deligne-Mumford stack associated to the
stacky fan $\mathbf{\Sigma}=(N,\Sigma,\beta)$. For a cone
$\sigma$ in the  simplicial fan $\Sigma$, let
$link(\sigma)=\{b_{i}: \rho_{i}+\sigma~ \text{is a cone in}~
\Sigma\}$. Then we have a quotient stacky fan
$\mathbf{\Sigma/\sigma}=(N(\sigma),\Sigma/\sigma,\beta(\sigma))$,
where
$$\beta(\sigma):\mathbb{Z}^{l}\to N(\sigma)$$
is given by the images of $\{b_{i}\}$'s in $link(\sigma)$.
Let $m:=|\sigma|$, then $dim(N_{\sigma})=m$ since $\sigma$ is simplicial.
Consider the commutative diagrams
\begin{equation}\label{close-hypertoric-a}
\begin{CD}
0 @ >>>\mathbb{Z}^{l+m}@ >>> \mathbb{Z}^{n}@ >>>
\mathbb{Z}^{n-l-m} @
>>> 0\\
&& @VV{\widetilde{\beta}}V@VV{\beta}V@VV{}V \\
0@ >>> N @ >{\cong}>>N@ >>> 0 @>>> 0,
\end{CD}
\end{equation}
and
\begin{equation}\label{close-hypertoric-aa}
\begin{CD}
0 @ >>>\mathbb{Z}^{m}@ >>> \mathbb{Z}^{l+m}@ >>>
\mathbb{Z}^{l} @
>>> 0\\
&& @VV{\beta_{\sigma}}V@VV{\widetilde{\beta}}V@VV{\beta(\sigma)}V \\
0@ >>> N_{\sigma} @ >{}>>N@ >>> N(\sigma) @>>> 0.
\end{CD}
\end{equation}
Applying the
Gale dual yields
\begin{equation}\label{close-hypertoric1}
\begin{CD}
0 @ >>>\mathbb{Z}^{n-l-m}@ >>> \mathbb{Z}^{n}@ >>>
\mathbb{Z}^{l+m} @
>>> 0\\
&& @VV{\cong}V@VV{\beta^{\vee}}V@VV{\widetilde{\beta}^{\vee}}V \\
0@ >>> \mathbb{Z}^{n-l-m} @ >{}>>N^{\vee}@ >{\phi_{1}}>> \widetilde{N}^{\vee}@>>>0,
\end{CD}
\end{equation}
and
\begin{equation}\label{close-hypertoric2}
\begin{CD}
0 @ >>>\mathbb{Z}^{l}@ >>> \mathbb{Z}^{l+m}@ >>>
\mathbb{Z}^{m} @
>>> 0\\
&& @VV{\beta(\sigma)^{\vee}}V@VV{\widetilde{\beta}^{\vee}}V@VV{\beta_{\sigma}^{\vee}}V \\
0@ >>> N(\sigma)^{\vee} @ >{\phi_{2}}>>\widetilde{N}^{\vee}@ >>> N_{\sigma}^{\vee}
@>>> 0.
\end{CD}
\end{equation}
Since $\mathbb{Z}^{m}\cong N_{\sigma}$, we have $N_{\sigma}^{\vee}=0$.
Applying $\text{Hom}_{\mathbb{Z}}(-,\mathbb{C}^{*})$ to (\ref{close-hypertoric1}), (\ref{close-hypertoric2}) yields
\begin{equation}\label{close-hypertoric3}
\begin{CD}
1 @ >>>\widetilde{G}@ >{\varphi_{1}}>> G@ >>>
(\mathbb{C}^{*})^{n-l-m} @
>>> 1\\
&& @VV{\widetilde{\alpha}}V@VV{\alpha}V@VV{\cong}V \\
1@ >>> (\mathbb{C}^{*})^{l+m} @ >{}>>(\mathbb{C}^{*})^{n}@ >>> (\mathbb{C}^{*})^{n-l-m}
@>>> 1,
\end{CD}
\end{equation}
and
\begin{equation}\label{close-toric4}
\begin{CD}
1 @ >>>1@ >>> \widetilde{G}@ >{\cong}>>
G(\sigma) @
>>> 1\\
&& @VV{}V@VV{\widetilde{\alpha}}V@VV{\alpha(\sigma)}V \\
1@ >>> (\mathbb{C}^{*})^{m} @ >{}>>(\mathbb{C}^{*})^{l+m}@ >>> (\mathbb{C}^{*})^{l}
@>>> 1.
\end{CD}
\end{equation}
Let $Z(\sigma)=\mathbb{A}^{l}\setminus \mathbb{V}(J_{\Sigma/\sigma})$, where
$J_{\Sigma/\sigma}$ is the irrelevant ideal of the quotient simplicial fan $\Sigma/\sigma$.
By the definition of toric Deligne-Mumford stack, we have
$\mathcal{X}(\mathbf{\Sigma/\sigma})=[Z(\sigma)/G(\sigma)]$,
where the action of $G(\sigma)$ is via the map $\alpha(\sigma)$
in (\ref{close-toric4}).

\begin{prop}\label{closedsub}
If $\sigma$ is a cone in the simplicial fan $\Sigma$,
then
$\mathcal{X}(\mathbf{\Sigma/\sigma})$ is a closed substack of
$\mathcal{X}(\mathbf{\Sigma})$.
\end{prop}

\begin{pf}
Let $W(\sigma)$ be the subvariety of $Z$ defined
by the ideal $J(\sigma):=\left<z_{i}:\rho_{i}\subseteq
\sigma\right>$. Then $W(\sigma)$ contains the $\mathbb{C}$-points
$z\in \mathbb{C}^{n}$ such that
the cone spanned by $\{\rho_{i}: z_{i}=0\}$ containing
$\sigma$ belongs to $\Sigma$. Then
the $\mathbb{C}$-point $z$ in $W(\sigma)$ such that
$\rho_{i}\nsubseteq \sigma\cup link(\sigma)$ implies that
$z_{i}\neq 0$. This implies that $W(\sigma)\cong Z(\sigma)\times (\mathbb{C}^{*})^{n-l-m}$. It is clear that $W(\sigma)$ is invariant under the $G$-action.

Let $\varphi_{0}: Z(\sigma)\to W(\sigma)$ be the inclusion
given by $z\mapsto (z,1)$. Then we have
a morphism of groupoids $\varphi_{0}\times \varphi_{1}: Z(\sigma)\times G(\sigma)\toto W(\sigma)\times G$ which
induces a morphism of stacks $\varphi: [Z(\sigma)/G(\sigma)]\to [W(\sigma)/G]$.
To prove that it is an isomorphism, we first prove that the following
diagram is cartesian:
$$\xymatrix{
Z(\sigma)\times G(\sigma)\dto_{(s,t)}\rto^{\varphi_{0}\times\varphi_{1}} &W(\sigma)\times G\dto^{(s,t)}\\
Z(\sigma)\times Z(\sigma)\rto^{\varphi_{0}\times\varphi_{0}}
&W(\sigma)\times W(\sigma).} $$
This is easy to prove. Given an element $(z_1,z_2)\in Z(\sigma)\times Z(\sigma)$, under
the map $\varphi_{0}\times\varphi_{0}$, we get $((z_1,1),(z_2,1))\in W(\sigma)\times W(\sigma)$.
If there is an element $g\in G$ such that $g(z_1,1)=(z_2,1)$, then from the exact sequence in the first row
of (\ref{close-hypertoric3}), there is an element $g(\sigma)\in G(\sigma)$ such that $g(\sigma)z_1=z_2$.
Thus we have an element $(z_1,g(\sigma))\in Z(\sigma)\times G(\sigma)$. So the
morphism $\varphi: [Z(\sigma)/G(\sigma)]\to [W(\sigma)/G]$ is injective.
Let $(z,1)$ be an element in $W(\sigma)$, then there exists an element $g\in (\mathbb{C}^{*})^{n-l-m}$ such that
$g(z,1)=(z,1)$. By (\ref{close-hypertoric3}), $g$ determines an element in $G$, so
$\varphi$ is surjective and $\varphi$ is an isomorphism.
Clearly  the stack $[W(\sigma)/G]$ is a closed substack of $\mathcal{X}(\mathbf{\Sigma})$, so $\mathcal{X}(\mathbf{\Sigma/\sigma})=[Z(\sigma)/G(\sigma)]$ is
also a closed substack of $\mathcal{X}(\mathbf{\Sigma})$.
\end{pf}

\begin{rmk}
Proposition \ref{closedsub} is Proposition 4.2 of \cite{BCS}. However the proof given there has a gap: some incorrect exact sequences were used. We choose to give a new proof here. One can prove this result using extended stacky fan defined in \cite{Jiang} since the quotient stacky fan
is naturally an extended stacky fan.
\end{rmk}

Following \cite{BCS}, for each top dimensional cone $\sigma$ in
$\Sigma$, denote by $Box(\sigma)$ the set of elements $v\in N$
such that $\overline{v}=\sum_{\rho_{i}\subseteq
\sigma}a_{i}\overline{b}_{i}$ for some $0\leq a_{i}<1$. The elements
in  $Box(\sigma)$ are in one-to-one correspondence with the elements
in the finite group $N(\sigma)=N/N_{\sigma}$, where $N(\sigma)$ is a
local  isotropy group of the stack $\mathcal{X}(\mathbf{\Sigma})$. If
$\tau\subseteq \sigma$ is a subcone, we define
$Box(\tau)$ to be the set of elements in $v\in N$ such that
$\overline{v}=\sum_{\rho_{i}\subseteq \tau}a_{i}\overline{b}_{i}$,
where $0\leq a_{i}<1$. It is easy to see that $Box(\tau)\subset
Box(\sigma)$. In fact the elements in $Box(\tau)$ generate a
subgroup of the local group $N(\sigma)$. Let
$Box(\mathbf{\Sigma})$ be the union of $Box(\sigma)$ for all
$d$-dimensional cones $\sigma\in \Sigma$. For $v_{1},\ldots,v_{n}\in
N$, let $\sigma(\overline{v}_{1},\ldots,\overline{v}_{n})$ be the
unique minimal cone in $\Sigma$ containing
$\overline{v}_{1},\ldots,\overline{v}_{n}$.

\begin{prop}[\cite{BCS}]\Label{r-inertia}
The $r$-th inertia stack of the stack $\mathcal{X}(\mathbf{\Sigma})$  is
$$\mathcal{I}_{r}\left(\mathcal{X}(\mathbf{\Sigma})\right)=\coprod_{(v_{1},\cdots,v_{r})\in Box(\mathbf{\Sigma})^{r}}
~\mathcal{X}(\mathbf{\Sigma/\sigma}(\overline{v}_{1},\cdots,\overline{v}_{r})).$$
\end{prop}

We are interested in the cases $r=1$ or $2$. When $r=1$,
$$\mathcal{I}\left(\mathcal{X}(\mathbf{\Sigma})\right)=\coprod_{v\in Box(\mathbf{\Sigma})}
~\mathcal{X}(\mathbf{\Sigma/\sigma}(\overline{v}))$$ is the
inertia stack. The orbifold Chow ring is the Chow ring of the
inertia stack as $\mathbb{Z}$-modules.

When $r=2$, for any pair $(v_{1},v_{2})$ in $Box(\mathbf{\Sigma})$, there is a unique
$v_{3}\in Box(\mathbf{\Sigma})$ such that $v_{1}+v_{2}+v_{3}\equiv 0\,\,(\text{mod }N)$. We have:
\begin{equation}\label{3-twisted sector}
\mathcal{I}_{2}\left(\mathcal{X}(\mathbf{\Sigma})\right)=\coprod_{(v_{1},v_{2},v_{3});v_{1}+v_{2}+v_{3}\equiv 0\,\,(\text{mod }N)}
~\mathcal{X}(\mathbf{\Sigma/\sigma}(\overline{v}_{1},\overline{v}_{2},\overline{v}_{3})).
\end{equation}
The components are called 3-twisted sectors in \cite{CR1}.

\subsection{The integral orbifold Chow ring.}

Let $\mathcal{X}(\mathbf{\Sigma})$ be a toric Deligne-Mumford stack with stacky fan
$\mathbf{\Sigma}$ and $A_{orb}^*(\mathcal{X}(\mathbf{\Sigma}),\mathbb{Z})$ its integral orbifold Chow ring. We first study the $A^*(\mathcal{X}(\mathbf{\Sigma}),\mathbb{Z})$-module structure of $A_{orb}^*(\mathcal{X}(\mathbf{\Sigma}),\mathbb{Z})$. Because $\Sigma$ is a simplicial fan, we have the following two lemmas in \cite{Jiang}:
\begin{lem}\label{lemma1}
For any $c\in N$, let $\sigma$ be the minimal cone in $\Sigma$
containing $\overline{c}$, then there exists a unique expression
$$c=v+\sum_{\rho_{i}\subset\sigma}m_{i}b_{i}$$
where $m_{i}\in \mathbb{Z}_{\geq 0}$, and $v\in Box(\sigma)$.
$\square$
\end{lem}

\begin{lem}\label{lemma2}
Let $\tau$ be a cone in the complete simplicial fan $\Sigma$ and
$\{\rho_{1},\ldots,\rho_{s}\}\subset link(\tau)$. Suppose
$\rho_{1},\ldots,\rho_{s}$ are contained in a cone $\sigma\subset
\Sigma$. Then $\sigma\cup\tau$ is contained in a cone of $\Sigma$.
$\square$
\end{lem}

Let $v\in Box(\mathbf{\Sigma})$ and
$\sigma:=\sigma(\overline{v})$ the minimal cone containing
$\overline{v}$. Then we have the
quotient stacky fan $\mathbf{\Sigma/\sigma}$ and
$\mathbf{\Sigma_{red}/\sigma}$. From the diagrams
(\ref{close-hypertoric-a}) and (\ref{close-hypertoric-aa}) we consider
the following diagrams:
\begin{equation}\label{firstdiag1}
\xymatrix@=5pt{
     & & ~\mathbb{Z}^{n} \ar[dr] \ar[ddd] &
      & & & & & & ~\mathbb{Z}^{l} \ar[dr] \ar[ddd] & \\
  ~\mathbb{Z}^{l+m} \ar[urr] \ar[dr] \ar[ddd] & & & N \ar[ddd]^{}
    & & & & ~\mathbb{Z}^{l+m} \ar[urr] \ar[dr] \ar[ddd] & & & N(\sigma) \ar[ddd]^{} \\
     & N \ar[urr]^{\cong} \ar[ddd] & & & \ar@{}[rr]^{}
       & & & & N \ar[urr] \ar[ddd] & & \\
     & & ~\mathbb{Z}^{n} \ar[dr] &
      & & & & & & ~\mathbb{Z}^{l} \ar[dr] & \\
  ~\mathbb{Z}^{l+m} \ar[urr] \ar[dr] & & & ~\overline{N}
   & & & & ~\mathbb{Z}^{l+m} \ar[urr] \ar[dr] & & & ~\overline{N}(\sigma). \\
     & ~\overline{N} \ar[urr]^{\cong} & &
     & & & & & ~\overline{N} \ar[urr] & &
  }
\end{equation}

Taking Gale dual yields

\begin{equation}\label{firstdiag2}
\xymatrix@=5pt{
     & & ~\mathbb{Z}^{l+m} \ar[dr] \ar[ddd] &
      & & & & & & ~\mathbb{Z}^{l+m} \ar[dr] \ar[ddd] & \\
  ~\mathbb{Z}^{n} \ar[urr] \ar[dr] \ar[ddd] & & & ~\widetilde{\overline{N}}^{\vee} \ar[ddd]^{\widetilde{\varphi}}
    & & & & ~\mathbb{Z}^{l} \ar[urr] \ar[dr] \ar[ddd] & & & ~\widetilde{\overline{N}}^{\vee} \ar[ddd]^{\widetilde{\varphi}} \\
     &~\overline{N}^{\vee} \ar[urr]^{} \ar[ddd]^{\varphi} & & & \ar@{}[rr]^{}
       & & & & ~\overline{N}(\sigma)^{\vee} \ar[urr] \ar[ddd]^{\varphi(\sigma)} & & \\
     & & ~\mathbb{Z}^{l+m} \ar[dr] &
      & & & & & & ~\mathbb{Z}^{l+m} \ar[dr] & \\
  ~\mathbb{Z}^{n} \ar[urr] \ar[dr] & & & ~\widetilde{N}^{\vee}
   & & & & ~\mathbb{Z}^{l} \ar[urr] \ar[dr] & & & ~\widetilde{N}^{\vee}. \\
     & ~N^{\vee} \ar[urr] & &
     & & & & & ~N(\sigma)^{\vee} \ar[urr] & &
  }
\end{equation}

For the quotient stacky fan $\mathbf{\Sigma/\sigma}$, if in the
map $\beta: \mathbb{Z}^{n}\to N$, the vectors
$\{b_1,\cdots,b_{n}\}$ generate the torsion part of $N$, then from 
(\ref{firstdiag1}) and (\ref{firstdiag2}),  
the vectors $\{\widetilde{b}_1,\cdots,\widetilde{b}_{l}\}$ in
the map $\beta(\sigma): \mathbb{Z}^{l}\to N(\sigma)$
generate the torsion part of $N(\sigma)$. So we have:

\begin{prop}\label{gerbequotient}
Given a toric Deligne-Mumford stack $\mathcal{X}(\mathbf{\Sigma})$
associated to the stacky fan $\mathbf{\Sigma}$. Suppose that the
map $\beta$ generate the torsion part of $N$, then for a cone
$\sigma\subset\Sigma$, the closed substack
$\mathcal{X}(\mathbf{\Sigma/\sigma})$ is a nontrivial gerbe over
the toric orbifold $\mathcal{X}(\mathbf{\Sigma_{red}/\sigma})$.
$\square$
\end{prop}

By Lemma \ref{diagonal}, the map $\varphi$ in (\ref{firstdiag2}) is diagonalizable,
so the map $\widetilde{\varphi}$ in (\ref{firstdiag2}) is diagonalizable.
From (\ref{close-hypertoric1}) and (\ref{close-hypertoric2}),
$N(\sigma)^{\vee}\cong \widetilde{N}^{\vee}$, and 
$\overline{N}(\sigma)^{\vee}\cong \widetilde{\overline{N}}^\vee$, so the 
map $\varphi(\sigma)$ in (\ref{firstdiag2}) is diagonalizable. We assume that 
$\varphi(\sigma)$ is given by the diagonal integral matrix $M(\sigma)$.
Let
$\widetilde{\mathbf{x}_{l}}=(\widetilde{x}_{1},\cdots,\widetilde{x}_{l})$
and $\mathbf{x}_{l}=(x_1,\cdots,x_l)$ be column vectors. Then using the same analysis
of the formula (\ref{key-formula}),  we have the following lemma:

\begin{lem}\label{lem-matrix}
The formula in (\ref{key-formula}) induces a formula
$$\widetilde{\mathbf{x}_{l}}=A(\sigma)M(\sigma)C(\sigma)\mathbf{x}_{l}$$
for the quotient stacky fan, where $A(\sigma)$ and $C(\sigma)$
are integral matrices. When we take $x_j$
as first Chern class of the line bundle $L_j$, 
the definition of $\tilde{x}_i$ are compatible with 
restrictions to components of the inertia stack. $\square$
\end{lem}

\begin{prop}\label{moduleisom}
Let $\mathcal{X}(\mathbf{\Sigma})$ be a
toric Deligne-Mumford stack associated to the stacky fan
$\mathbf{\Sigma}$, then we have an isomorphism of
$A^{*}(\mathcal{X}(\mathbf{\Sigma}),\mathbb{Z})$-modules:
$$\bigoplus_{v\in Box(\mathbf{\Sigma})}A^{*}\left(\mathcal{X}(\mathbf{\Sigma/\sigma}(\overline{v})),\mathbb{Z}\right)[deg(y^{v})]\cong
\frac{\mathbb{Z}[\mathbf{\Sigma}]}{Cir(\mathbf{\Sigma})}.$$
\end{prop}
\begin{pf}
We use a method similar to that in Proposition 4.7 of \cite{Jiang}. Let
$$S_{\mathbf{\Sigma}}:=\mathbb{Z}[y^{b_{i}}:
\rho_{i}\in\Sigma(1)]/\mathcal{I}_{\mathbf{\Sigma}}.$$ Then
$S_{\mathbf{\Sigma}}/Cir(\mathbf{\Sigma})\cong
A^{*}(\mathcal{X}(\mathbf{\Sigma}),\mathbb{Z})$ given by
$y^{b_{i}}\mapsto x_i$. By the definition of $\mathbb{Z}[\mathbf{\Sigma}]$ and Lemma \ref{lemma1}, we see  that $\mathbb{Z}[\mathbf{\Sigma}]=\bigoplus_{v\in
Box(\mathbf{\Sigma})}y^{v}\cdot S_{\mathbf{\Sigma}}$. And we obtain an
isomorphism of
$A^{*}(\mathcal{X}(\mathbf{\Sigma}),\mathbb{Z})$-modules:
\begin{equation}\Label{module}
\frac{\mathbb{Z}[\mathbf{\Sigma}]}{Cir(\mathbf{\Sigma})}\cong
\bigoplus_{v\in Box(\mathbf{\Sigma})}\frac{y^{v}\cdot
S_{\mathbf{\Sigma}}}{y^{v}\cdot Cir(\mathbf{\Sigma})}.
\end{equation}

For any  $v\in Box(\mathbf{\Sigma})$,  let $\sigma(\overline{v})$ be
the minimal cone in $\Sigma$ containing $\overline{v}$.  Let
$\rho_{1},\ldots, \rho_{l}\in link(\sigma(\overline{v}))$, and
$\widetilde{\rho}_{i}$ be the image of $\rho_{i}$ under the natural
map $\pi: N\to
N(\sigma(\overline{v}))=N/N_{\sigma(\overline{v})}$. Then
$S_{\mathbf{\Sigma}/\sigma(\overline{v})}\subset
\mathbb{Z}[\mathbf{\Sigma}/\sigma(\overline{v})]$
is the subring generated by: $y^{\widetilde{b}_{i}}$, for
$\rho_{i}\in link(\sigma(\overline{v}))$. Let $\widetilde{a}$ be the
order of the torsion subgroup of $N(\sigma(\overline{v}))$. Then let
$a=s\widetilde{a}$, and conversely we have
$\widetilde{a}=\frac{1}{s}a$. By Lemmas \ref{lemma2} and \ref{lem-matrix}, it is easy to check that the ideal 
$\mathcal{I}_{\mathbf{\Sigma}/\sigma(\overline{v})}$ goes to the ideal $\mathcal{I}_{\mathbf{\Sigma}}$ and we have a morphism $\widetilde{\Psi}_{v}:
S_{\mathbf{\Sigma}/\sigma(\overline{v})}[deg(y^{v})]\to
y^{v}\cdot S_{\Sigma}$  given by: $y^{\widetilde{b}_{i}}\mapsto
y^{v}\cdot sy^{b_{i}}$. If
$\sum_{i=1}^{l}\widetilde{\theta}(\widetilde{b}_{i})\widetilde{a}y^{\widetilde{b}_{i}}$
belongs to the ideal $Cir(\mathbf{\Sigma/\sigma}(\overline{v}))$,
then
$$
\widetilde{\Psi}_{v}\left(\sum_{i=1}^{l}\widetilde{\theta}(\widetilde{b}_{i})\widetilde{a}y^{\widetilde{b}_{i}}\right)
= y^{v}\cdot
\left(\sum_{i=1}^{n}\theta(b_{i})s\widetilde{a}y^{b_{i}}\right)
= y^{v}\cdot
\left(\sum_{i=1}^{n}\theta(b_{i})ay^{b_{i}}\right),
$$
where $\theta=\widetilde{\theta}\circ\pi$
and $\theta(b_{i})=\widetilde{\theta}(\widetilde{b}_{i})$. So we obtain that
$\widetilde{\Psi}_{v}(\sum_{i=1}^{l}\widetilde{\theta}(\widetilde{b}_{i})\widetilde{a}y^{\widetilde{b}_{i}})\in
y^{v}\cdot Cir(\mathbf{\Sigma})$. So
$\widetilde{\Psi}_{v}$ induce  a morphism $\Psi_{v}:
\frac{S_{\mathbf{\Sigma}/\sigma(\overline{v})}}{Cir(\mathbf{\Sigma}/\sigma(\overline{v}))}[deg(y^{v})]
\to \frac{y^{v}\cdot
S_{\mathbf{\Sigma}}}{y^{v}\cdot
Cir(\mathbf{\Sigma})}
$
such that $\Psi_{v}([y^{\widetilde{b}_{i}}])=[y^{v}\cdot
sy^{b_{i}}]$.

Conversely, for such $v\in Box(\mathbf{\Sigma})$ and
$\rho_{i}\subset \sigma(\overline{v})$, choose $\theta_{i}\in
\text{Hom}_\mathbb{Z}(N,\mathbb{Q})$ such that $\theta_{i}(b_{i})=1$ and
$\theta_{i}(b_{i^{'}})=0$ for $b_{i^{'}}\neq b_{i}\in
\sigma(\overline{v})$. Again by Lemmas \ref{lemma2} and \ref{lem-matrix} we consider the following morphism
$\widetilde{\Phi}_{v}: y^{v}\cdot S_{\mathbf{\Sigma}}\to
S_{\Sigma/\sigma(\overline{v})}[deg(y^{v})]$ given by:
$$y^{b_{i}}\mapsto\begin{cases}\frac{1}{s}y^{\widetilde{b}_{i}}&\text{if $\rho_{i}\subseteq link(\sigma(\overline{v}))$}\,,\\
-\sum_{j=1}^{l}\theta_{i}(b_{j})y^{\widetilde{b}_{j}}&\text{if
$\rho_{i}\subseteq \sigma(\overline{v})$}\,,\\
0&\text{if $\rho_{i}\nsubseteq \sigma(\overline{v})\cup
link(\sigma(\overline{v}))$}\,.\end{cases}$$
Let $y^{v}\cdot \left(\sum_{i=1}^{n}\theta(b_{i})ay^{b_{i}}\right)$
belong to the ideal $y^{v}\cdot Cir(\mathbf{\Sigma})$.
For $\theta\in M$, we have $\theta=\theta_{v}+\theta_{v}^{'}$,
where $\theta_{v}\in N(\sigma(\overline{v}))^{*}=M\cap
\sigma(\overline{v})^{\perp}$ and  $\theta_{v}^{'}$ belongs to the
orthogonal complement of the subspace
$\sigma(\overline{v})^{\perp}$ in $M$.  We have
$$
\widetilde{\Phi}_{v}\left(y^{v}\cdot
\left(\sum_{i=1}^{n}\theta(b_{i})ay^{b_{i}}\right)\right)
=\sum_{i=1}^{l}\theta_{v}(\widetilde{b}_{i})\widetilde{a}y^{\widetilde{b}_{i}}+
\sum_{\rho_{i}\subset
\sigma(\overline{v})}\theta_{v}^{'}(b_{i})\left(-\sum_{j=1}^{l}\theta_{i}(b_{j})y^{\widetilde{b}_{j}}\right)
+
\sum_{i=1}^{l}\theta_{v}^{'}(b_{i})y^{\widetilde{b}_{i}}.
$$
Note that
$\sum_{i=1}^{l}\theta_{v}(\widetilde{b}_{i})\widetilde{a}y^{\widetilde{b}_{i}}\in
Cir(\mathbf{\Sigma}/\sigma(\overline{v}))$. Now let
$\theta_{v}^{'}=\sum_{\rho_{i}\subset
\sigma(\overline{v})}a_{i}\theta_{i}$, where $a_{i}\in
\mathbb{Q}$,  then $\sum_{\rho_{i}\subset
\sigma(\overline{v})}\theta_{v}^{'}(b_{i})=\sum_{\rho_{i}\subset
\sigma(\overline{v})}a_{i}\theta_{i}(b_{i})$. We have:
$$\sum_{\rho_{i}\subset
\sigma(\overline{v})}a_{i}\theta_{i}(b_{i})\left(-\sum_{j=1}^{l}\theta_{i}(b_{j})y^{\widetilde{b}_{j}}\right)
+\sum_{\rho_{i}\subset
\sigma(\overline{v})}\sum_{j=1}^{l}a_{i}\theta_{i}(b_{j})y^{\widetilde{b}_{j}}=0,$$
so we have $\widetilde{\Phi}_{v}\left(y^{v}\cdot
\left(\sum_{i=1}^{n}\theta(b_{i})ay^{b_{i}}\right)\right)\in
Cir(\mathbf{\Sigma}/\sigma(\overline{v}))$. So
$\widetilde{\Phi}_{v}$ induces a morphism
$$\Phi: \frac{y^{v}\cdot
S_{\mathbf{\Sigma}}}{y^{v}\cdot
Cir(\mathbf{\Sigma})}\to
\frac{S_{\mathbf{\Sigma}/\sigma(\overline{v})}}{Cir(\mathbf{\Sigma}/\sigma(\overline{v}))}[deg(y^{v})].$$
We check that  $\Phi_{v}\Psi_{v}=1$ and $\Psi_{v}\Phi_{v}=1$. So $\Phi_{v}$ is
an isomorphism.
Note that both sides of (\ref{module}) are
$S_{\mathbf{\Sigma}}\slash Cir(\mathbf{\Sigma})=A^{*}(\mathcal{X}(\mathbf{\Sigma}),\mathbb{Z})$-modules,
we complete the proof.
\end{pf}

Now we compute the ring structure. The key part of the orbifold
cup product is the orbifold obstruction bundle. For the toric Deligne-Mumford stack
$\mathcal{X}(\mathbf{\Sigma})$, the obstruction bundle over the
3-twisted sectors in (\ref{3-twisted sector}) is given by:
\begin{prop}\label{obstr}
Let $\mathcal{X}(\mathbf{\Sigma/\sigma}(\overline{v}_{1},\overline{v}_{2},\overline{v}_{3}))$
be a 3-twisted sector of the toric Deligne-Mumford stack
$\mathcal{X}(\mathbf{\Sigma})$. Let
$v_{1}+v_{2}+v_{3}=\sum_{\rho_{i}\subset
\sigma(\overline{v}_{1},\overline{v}_{2},\overline{v}_{3})}a_{i}b_{i}$,
$a_{i}=1 , 2$, then the Euler class of the obstruction bundle
$O_{(v_{1},v_{2},v_{3})}$  over
$\mathcal{X}(\mathbf{\Sigma})_{(v_{1},v_{2},v_{3})}$ is:
$$\prod_{a_{i}=2}c_{1}(\mathcal{L}_{i})|_{\mathcal{X}(\mathbf{\Sigma/\sigma}(\overline{v}_{1},\overline{v}_{2},\overline{v}_{3}))},$$
where $\mathcal{L}_{i}$ is the line bundle over
$\mathcal{X}(\mathbf{\Sigma})$ corresponding to the ray $\rho_{i}$.
\end{prop}

\subsection*{Proof of Theorem \ref{chowring-main2}.}
Proposition \ref{moduleisom} gives an isomorphism between
$A^{*}(\mathcal{X}(\mathbf{\Sigma}), \mathbb{Z})$-modules:
$$A^{*}_{orb}\left(\mathcal{X}(\mathbf{\Sigma}),\mathbb{Z}\right)\simeq\bigoplus_{v\in
Box(\mathbf{\Sigma})}A^{*}\left(\mathcal{X}(\mathbf{\Sigma/\sigma}(\overline{v})),\mathbb{Z}\right)[deg(y^{v})]\cong
\frac{\mathbb{Z}[\mathbf{\Sigma}]}{Cir(\mathbf{\Sigma})}.$$ All we
need is  to show that the orbifold cup product defined in
\cite{AGV} coincides with the product in ring
$\mathbb{Z}[\mathbf{\Sigma}]/Cir(\mathbf{\Sigma})$. From the above
isomorphisms, it suffices to consider the canonical generators
$y^{b_{i}}$, $y^{v}$ where $v\in Box(\mathbf{\Sigma})$.

For any $v_{1},v_{2}\in Box(\mathbf{\Sigma})$, let $v_{3}\in Box(\mathbf{\Sigma})$
be the unique box element such that $v_{1}+v_{2}+v_{3}\equiv 0\,\,(\text{mod }N)$. Then
$\mathcal{X}(\mathbf{\Sigma/\sigma}(\overline{v}_{1},\overline{v}_{2},\overline{v}_{3}))$
is a 3-twisted sector. Let $e_{i}: \mathcal{X}(\mathbf{\Sigma/\sigma}(\overline{v}_{1},\overline{v}_{2},\overline{v}_{3}))
\to \mathcal{X}(\mathbf{\Sigma/\sigma}(\overline{v}_{i}))$ be the evaluation map
for $1\leq i\leq 3$. Let $\check{v}_{3}$ be the inverse of $v_{3}$
in the local group, and $I:\mathcal{X}(\mathbf{\Sigma/\sigma}(\overline{v}_{3}))\to
\mathcal{X}(\mathbf{\Sigma/\sigma}(\overline{\check{v}}_{3}))$ be the map given by
$(x,v_{3})\longmapsto (x,\check{v}_{3})$. Let $\check{e}_{3}=I\circ e_{3}$. Then the orbifold cup product is defined by:
$$y^{v_{1}}\cup_{orb}y^{v_{2}}=\check{e}_{3,*}(e^{*}_{1}y^{v_{1}}\cup e^{*}_{2}y^{v_{2}}\cup e(O_{(v_{1},v_{2},v_{3})})),$$
where $O_{(v_{1},v_{2},v_{3})}$ is the obstruction bundle in
Proposition \ref{obstr}. Since $e_{1},e_{2},\check{e}_{3}$ are all
inclusion, so are representable as morphisms of Deligne-Mumford
stacks. By \cite{Kr}, the pullback and pushforward are
well-defined for integral Chow classes. Let $\pi:
\mathcal{X}(\mathbf{\Sigma})\to\mathcal{X}(\mathbf{\Sigma_{red}})$
be the natural morphism of rigidification. The first Chern class
of the line bundle $L_{i}$ is $y^{b_{i}}$, so by Proposition
\ref{linebundle-com}, the first Chern class of $\mathcal{L}_{i}$
is $\sum_{j=1}^{n}a_{j,i}y^{b_{i}}$ which is
$\widetilde{y}^{b_{i}}$. This class  represents an integral Chow
class of
$\mathcal{X}(\mathbf{\Sigma/\sigma}(\overline{v}_{1},\overline{v}_{2},\overline{v}_{3}))$.
We have that $\check{e}_{3,*}(y^{b_{i}})=y^{b_{i}}\in
A^{*}(\mathcal{X}(\mathbf{\Sigma/\sigma}(\overline{\check{v}}_{3})),\mathbb{Z})$.
So by the definition of orbifold cup product we have
$$y^{v_{1}}\cdot y^{v_{2}}=
y^{\check{v}_{3}}\prod_{i\in I}\widetilde{y}^{b_{i}}\cdot
\prod_{i\in J}\widetilde{y}^{b_{i}}.$$ $\square$

\section{Examples}\label{example}
In this section we compute some examples of the integral Chow ring and integral orbifold Chow
rings.
\begin{example}[The moduli stack of 1-pointed elliptic curves]
Let $\Sigma$ be the complete fan of the projective line,
$N=\mathbb{Z}\oplus \mathbb{Z}/2\mathbb{Z}$, and $\beta:
\mathbb{Z}^{2}\to \mathbb{Z}\oplus
\mathbb{Z}/2\mathbb{Z}$ be given by the vectors
$\{b_{1}=(2,1),b_{2}=(-3,0)\}$. Then
$\mathbf{\Sigma}=(N,\Sigma,\beta)$ is a stacky fan.  We compute
that $(\beta)^{\vee}: \mathbb{Z}^{2}\to
N^{\vee}=\mathbb{Z}$ is given by the matrix [6,4]. So we get the
following exact sequence:
$$
0\longrightarrow \mathbb{Z}\longrightarrow
\mathbb{Z}^{2}\stackrel{\beta}{\longrightarrow}
\mathbb{Z}^{2}\oplus\mathbb{Z}_{2}\longrightarrow
0, $$
$$
0\longrightarrow \mathbb{Z}\longrightarrow
\mathbb{Z}^{2}\stackrel{\beta^{\vee}}{\longrightarrow}
\mathbb{Z}\longrightarrow\mathbb{Z}_{2}\longrightarrow
0, $$
and
\begin{equation}\Label{gerbe1}
1\longrightarrow \mu_{2}\longrightarrow
\mathbb{C}^{*}\stackrel{[6,4]^{t}}{\longrightarrow}
(\mathbb{C}^{*})^{2}\longrightarrow
\mathbb{C}^{*}\longrightarrow 1. \
\end{equation}
The toric Deligne-Mumford stack
$\mathcal{X}(\mathbf{\Sigma})=[\mathbb{C}^{2}-\{0\}/\mathbb{C}^{*}]=:\mathbb{P}(6,4)$,
where the action is given by
$\lambda(x,y)=(\lambda^{6}x,\lambda^{4}y)$, may be identified with the moduli stack
$\overline{\mathcal{M}}_{1,1}$ of 1-pointed elliptic curves. The
stack $\mathcal{X}(\mathbf{\Sigma})$ is the nontrivial
$\mu_{2}$-gerbe over $\mathbb{P}(3,2)$ coming from the canonical
line bundle over $\mathbb{P}(3,2)$. Since
$N=\mathbb{Z}\oplus\mathbb{Z}_{2}$, we have $m_{i}=2$. By Theorem
1.1, we have
$$A^{*}(\mathcal{X}(\mathbf{\Sigma}),\mathbb{Z})\cong \frac{\mathbb{Z}[x_{1},x_{2}]}
{(2x_{1}-3x_{2},2x_{1}2x_{2})}\cong\mathbb{Z}[t]/(24t^{2}),$$ which
is the same as the result in \cite{EG}. We compute the integral
orbifold Chow ring. There are 7 box elements:
$v=(1,1),w_{1}=(-1,0),w_{2}=(-2,0)$,$u=(0,1)$,$\rho_1=(1,0)$,$\rho_2=(-1,1)$ and $\rho_3=(-2,1)$ corresponding to 7
twisted sectors. The three box elements $v,w_1,u$ generate the others. So by Theorem 1.2 we have
\begin{eqnarray}
A_{orb}^{*}(\mathcal{X}(\mathbf{\Sigma}),\mathbb{Z})&\cong& \frac{\mathbb{Z}[x_{1},x_{2},v,w_{1},u]}
{(2x_{1}-3x_{2},2x_{1}2x_{2},v^{2}-2x_{1}u,w_{1}^{3}-2x_{2}u,vw_{1},v2x_{2},w_{1}2x_{1},u^{2}-1)} \nonumber \\
&\cong &\frac{\mathbb{Z}[t,v,w_{1},u]}
{(24t^{2},v^{2}-6tu,w_{1}^{3}-4tu,vw_{1},4vt,6w_{1}t,u^{2}-1)}\nonumber,
\end{eqnarray}
which is the same as the result in \cite{AGV}.
\end{example}

\begin{example}
In this example we discuss the relation between integral orbifold Chow ring and the integral Chow ring of crepant resolutions.
Let $N=\mathbb{Z}^{2}$, and $\beta:
\mathbb{Z}^{3}\to \mathbb{Z}^{2}$ be given by the vectors
$\{b_{1}=(1,0),b_{2}=(0,1),b_{3}=(-1,-2)\}$. Let $\Sigma$ be the complete fan in
$\mathbb{R}^{2}$ generated by $\{b_{1},b_{2},b_{3}\}$. Then
$\mathbf{\Sigma}=(N,\Sigma,\beta)$ is a stacky fan.  We compute
that $(\beta)^{\vee}: \mathbb{Z}^{3}\to
N^{\vee}=\mathbb{Z}$ is given by the matrix [1,1,2]. So we get the
following exact sequence:
\begin{equation}\Label{gerbe2}
1\longrightarrow
\mathbb{C}^{*}\stackrel{[1,1,2]^{t}}{\longrightarrow}
(\mathbb{C}^{*})^{3}\longrightarrow
(\mathbb{C}^{*})^{2}\longrightarrow 1 \
\end{equation}
The toric Deligne-Mumford stack
$\mathcal{X}(\mathbf{\Sigma})=[\mathbb{C}^{3}-\{0\}/\mathbb{C}^{*}]$,
where the action is given by $\lambda(x,y,z)=(\lambda x,\lambda
y,\lambda^{2}z)$, is the
weighted projective stack $\mathbb{P}(1,1,2)$ which is a toric
orbifold.  We compute the integral orbifold Chow ring. There is
one box element: $v=\frac{1}{2}b_{1}+\frac{1}{2}b_{3}=(0,1)$
corresponding to one twisted sector. So by Theorem 1.2 we have
\begin{eqnarray}
A_{orb}^{*}(\mathcal{X}(\mathbf{\Sigma}),\mathbb{Z})&\cong& \frac{\mathbb{Z}[x_{1},x_{2},x_{3},v]}
{(x_{1}x_{2}x_{3},x_{1}-x_{3},x_{2}-2x_{3},v^{2}-x_{1}x_{3},vx_{1},vx_{2},vx_{3})} \nonumber \\
&\cong &\frac{\mathbb{Z}[x_{3},v]}
{(2x_{3}^{3},2vx_{3},v^{2}-x_{3}^{2})}\nonumber.
\end{eqnarray}

Let $\rho_{4}$ be a ray generated by $v=b_{4}$. Then the complete fan $\Sigma^{'}=\{b_{1},b_{2},b_{3},b_{4}\}$
generated by the rays $\{\rho_{1},\rho_{2},\rho_{3},\rho_{4}\}$ is the fan of Hirzebruch surface
$\mathbb{F}_{2}$. It is well-known that 
\begin{eqnarray}
A^{*}(\mathbb{F}_{2},\mathbb{Z})&\cong& \frac{\mathbb{Z}[x_{1},x_{2},x_{3},x_{4}]}
{(x_{1}x_{2}x_{3},x_{1}-x_{3},x_{2}-2x_{3}-x_{4},,x_{2}x_{4},x_{1}x_{3})} \nonumber \\
&\cong &\frac{\mathbb{Z}[x_{3},x_{4}]}
{(2x_{3}^{3}+x_{3}^{2}x_{4},x_{3}^{2},2x_{3}x_{4}+x_{4}^{2})}\nonumber.
\end{eqnarray}
It is easy to see that these two rings are not isomorphic. So under integer coefficients,
$A_{orb}^{*}(\mathcal{X}(\mathbf{\Sigma}),\mathbb{Z})\ncong A^{*}(\mathbb{F}_{2},\mathbb{Z})$.
In \cite{BMP}, the authors proved that $$A_{orb}^{*}(\mathcal{X}(\mathbf{\Sigma}),\mathbb{C})
\cong A^{*}_{Q}(\mathbb{F}_{2},\mathbb{C})$$
where $A^{*}_{Q}(\mathbb{F}_{2},\mathbb{Z})$ is the quantum corrected cohomology
of $\mathbb{F}_{2}$ under complex number coefficients, thus verifying the {\em cohomological crepant resolution conjecture} \cite{R}. 
\end{example}


\end{document}